\documentclass[reqno]{amsart}

\theoremstyle{definition}

\theoremstyle{remark}

\numberwithin{equation}{section}

\usepackage[all]{xy}

\newtheorem{tm}{Theorem}[subsection]
\newtheorem{lm}[tm]{Lemma}
\newtheorem{pr}[tm]{Proposition}
\newtheorem{rmk}[tm]{Remark}

\newtheorem{ex}[tm]{Example}

\newtheorem{??}[tm]{Question}
\newtheorem{defi}[tm]{Definition}

\font\tenmsb=msbm10
\font\sevenmsb=msbm7
\font\fivemsb=msbm5
    
\newfam\msbfam
\textfont\msbfam=\tenmsb
\scriptfont\msbfam=\sevenmsb
\scriptscriptfont\msbfam=\fivemsb
\def\Bbb#1{{\fam\msbfam #1}}

\font\teneufm=eufm10
\font\seveneufm=eufm7
\font\fiveeufm=eufm5
\newfam\eufmfam
\textfont\eufmfam=\teneufm
\scriptfont\eufmfam=\seveneufm
\scriptscriptfont\eufmfam=\fiveeufm
\def\frak#1{{\fam\eufmfam\relax#1}}

\def\lorw{\longrightarrow}
\newcommand\n{\noindent}

\newcommand\ci{\cite}

\newcommand\rat{{\Bbb Q}}

\newcommand\zed{{\Bbb Z}}

\newcommand\blacksquare{{\hspace*{\fill} $\fbox{}$}}

\newcommand{\im}{ \hbox{\rm Im} }

\newcommand{\ke}{ \hbox{\rm Ker} }

\newcommand{\ptd}[1]{ \,^{\frak p}\!\tau_{ \leq {#1} } }

\newcommand{\ptu}[1]{ \,^{\frak p}\!\tau_{ \geq {#1} } }

\newcommand{\td}[1]{ \tau_{ \leq {#1} } }

\newcommand{\pp}{ \,^{\frak p}\!\tau    }
\newcommand{\ts}[1]{ \tau^{ #1  }}
\newcommand{\pe}{ {\mathcal P }  }
\newcommand{\tb}[1]{ \tau_{ #1  }}
\newcommand{\pc}[2]{ \,^{\frak p}\!{\mathcal H}^{#1}({#2})   }
\newcommand{\tu}[1]{ \tau_{ \geq {#1} } }

\begin{document}

\title[The perverse filtration and the Lefschetz hyperplane theorem, II]{The perverse filtration  and \\  the  Lefschetz hyperplane theorem, II}
\author{
Mark Andrea A.  de Cataldo
}
\address{Department of Mathematics, Stony Brook University, Stony Brook, NY 11794}
\email{mde@math.sunysb.edu}
\thanks{The author was supported in part by  NSA and NSF grants.}



\dedicatory{This paper is dedicated to Prakob Monkolchayut.}


\begin{abstract}
The perverse filtration in cohomology and in cohomology with compact supports
is interpreted in terms of kernels of restrictions maps
to suitable subvarieties by  using the  Lefschetz hyperplane theorem and  spectral objects.
Various mixed-Hodge-theoretic consequences  for intersection cohomology
and for the decomposition theorem are derived.
\end{abstract}

\maketitle

\tableofcontents

\section{Introduction}
\label{secintro}
Let $Y$ be an affine variety of dimension $n$, let $K$ be a bounded complex
of sheaves of abelian groups on $Y$ with constructible cohomology sheaves,
 and let $H^*(Y, K)$ be the (hyper)cohomology groups of $Y$ with coefficients in $K$.
 For simplicity, in this introduction, we confine ourselves
to affine varieties and to cohomology. In this paper,   we prove analogous  results for quasi projective varieties,
  and  for cohomology 
with compact supports. 
Fix an arbitrary embedding $Y \subseteq {\Bbb A}^N$ in affine space and let  
 $Y_* = \{ Y_{-n} \subseteq \ldots \subseteq Y_0 = Y \}$ be a  general $n$-flag of linear sections  of $Y$,
 i.e. $Y_{-i} = Y \cap \Lambda^{N-i}$, where  $\{\Lambda^{N-n} \subseteq \ldots  \subseteq
 \Lambda^{N} = {\Bbb A}^N\}$ is a general 
 partial flag of linear affine subspaces of ${\Bbb A}^N$.

In  \ci{decmigpf}, we showed
that the (middle)  perverse filtration on the cohomology 
groups $H^*(Y,K)$  (\S\ref{ssawts})
can be described
geometrically as follows:  up to renumbering,  {\em the   perverse filtration  coincides
with the  flag filtration}, i.e. the subspaces of the (middle) perverse filtration
coincide with the
 kernels of  the restriction maps $H^*(Y,K) \to H^*(Y_p, K_{|Y_p})$.

The purpose of this paper is twofold.

The  former is to use Verdier's spectral objects
to  give, in $\S$\ref{tgofsss},  an alternative proof   of  the main result
of \ci{decmigpf}, i.e. of the description of the perverse filtration in cohomology
and in cohomology with compact supports using 
general flags. 
In order to do so, I introduce a new technique  that,
in the presence of suitable $t$-exactness, 
  ``realigns"  a map of spectral sequences; see Lemma \ref{sotex}.
I hope 
   that this technique is of independent  interest and 
  will have  further applications.

The latter is to use classical  mixed  Hodge theory
and  the description of the perverse filtration via flags
to establish in a rather elementary
fashion, in $\S$\ref{dtamhs},
a rather complete mixed-Hodge-theoretic package for the  intersection cohomology groups
of quasi projective varieties and for the maps between them.
  For a survey on the
decomposition theorem, see \ci{decmigbams}.
The results of this paper construct the  mixed Hodge structures in a direct way,
by using  Deligne's theory, the decomposition theorem and various geometric constructions.

M. Saito proved these results in \ci{samhm} by  using mixed Hodge modules.
While the   mixed Hodge structures in $\S$\ref{dtamhs}  could  be a priori  different from the ones
stemming from M. Saito's work,  Theorem \ref{coincide} asserts that
the two structures coincide.

The mixed-Hodge theoretic results are stated in $\S$\ref{dtamhs} and they are proved
in $\S$\ref{schema} and $\S$\ref{-1}. The strategy is to first prove
the results in the case when the domain $X$  of the map $f: X \to Y$ is 
nonsingular, and then to  use the nonsingular case and resolution of singularities.
This strategy is adapted from \ci{decmightam, decmigseattle}, which deal with the
case of projective varieties. In the projective case, the intersection cohomology groups
coincide with the ones with compact supports,  the Hodge structures are pure
and  key ingredients are   the use of the intersection pairing on intersection cohomology
and a {\em different} geometric description of the perverse filtration (\ci{decmightam}), valid {\em
only} in the projective case.
In the present paper, 
we make a systematic use of the description of the perverse filtration
via flags and, since we deal with non compact varieties, the intersection pairing
involves compact supports as well, and
we must deal with  intersection cohomology and with intersection cohomology
with compact supports simultaneously. 
In the course of the proofs, care has to be taken to verify many mixed-Hodge-theoretic compatibilities, and
I have included  the details of the proof that I feel are not entirely a matter of routine.

Let me  try to describe the meaning and the usefulness of
the ``realignment'' Lemma \ref{sotex} by discussing a simple situation; see also Remark \ref{triv}.

Let   $P$ be a perverse sheaf on  $Y={\Bbb A}^n$
 and  $i: Z\to Y$ be a general hyperplane. 
  In general,  while $P_{|Z}$ is not a perverse sheaf on $Z$,
the shifted $P[-1]_{|Z}$ is.
 The natural adjunction  map
$P \to i_*i^* P$ yields the restriction map
$r: H^*(Y, P) \to H^*(Z, P_{|Z})$. 
Since  the cohomology perverse sheaves $\pc{t}{P} =0$ for $t\neq 0$ and $\pc{t}{P_{|Z}} =0$ for $t\neq -1$,
the  map of perverse  spectral sequences
$r_2^{st}: H^s(Y, \pc{t}{P}) \to H^s (Z, \pc{t}{P_{|Z}})$
that arises naturally by functoriality  is the zero map.
On the other hand, the Lefschetz hyperplane theorem for perverse sheaves
implies that the restriction map  $r$ is injective in negative cohomological 
degrees. 

The conclusion is that the map of spectral sequences
that arises from restriction from $Y$ to $Z$ does not adequately reflect
the geometry. In technical
   jargon, the issue  is that 
  $i^*$ is not $t$-exact and the map of spectral sequences above
  has  components  $r_2^{st}$ relating the ``wrong" groups.  
This does not happen
  in the standard case, for $i^*$ is always exact; see Remark \ref{flss}.
  
The statement of the Lefschetz hyperplane theorem for perverse sheaves, 
coupled with the identities $P= \pc{0}{P}$ and  $P_{|Z} = \pc{-1}{P_{|Z}}[1]$,
suggests  that we should try and
tilt  the arrows $r_2^{st}$ so that they start  at a spot $(s,t)$ but
end at a spot $(s+1, t-1)$, thus being the arrows appearing in the 
Lefschetz hyperplane theorem. This tilting of the arrows can also  
be seen as a realignment of the pages of the target perverse spectral sequence
(for $Z$)  and we
say that the target spectral sequence has been {\em translated} by one unit
(see (\ref{tsspe})).

The point is that this should be achieved in 
a coherent way, i.e.  we are looking for a new map of spectral sequences,
which has as source the original source, and as target the translate by one unit of the original target.
This is precisely what Lemma \ref{sotex} achieves
in the more general context of spectral objects in $t$-categories.
In the special case discussed above, the key point is that
 $i^*[-1]$ is suitably $t$-exact in view of the fact that the hyperplane $H$ is general
 (see Remark \ref{1122}),
 hence transversal to all the strata of a stratification for $P$.
  This explains
 why   we need to choose the flag of linear sections $Y_*$ to be general.
The perverse spectral sequence
for $H^*(Y,P)$ is now mapped nontrivially to the $(-1)$-translate of the perverse spectral sequence
for $H^*(Z, P_{|Z})$. 

The  geometric description of the perverse filtration using
flags  (Theorem \ref{mtm})   is based on an iteration  
of  this procedure for each element of the flag.
By doing so, we obtain a compatible system of realigned maps of spectral sequences
whose  first pages have many zeroes (Artin vanishing theorem \ref{cdav}) 
and are connected by (mostly) monic arrows
-epic if we are dealing with  compact supports- (Lefschetz hyperplane theorem
\ref{swl}). The identification of the flag and perverse filtrations is then carried out by inspecting
the pages of the realigned maps of spectral sequences (see $\S$\ref{alossq} and 
the proof of Theorem \ref{mtm}).

\bigskip
{\bf Acknowledgments.}
It is a pleasure to thank D. Arapura,  A. Beilinson, M. Goresky, 
M. Levine, M. Nori 
for stimulating conversations, and the anonymous referee for the  very careful work.
Special thanks go to L. Migliorini:
the results of $\S$\ref{secappl} have been obtained jointly with him.
I thank
the University of Bologna, the I.A.S., Princeton, and Max Planck
Institut in Bonn
for their hospitality during the preparation of parts of this paper.

\section{Notation and background results}
\label{notandback}
This paper is a sequel of \ci{decmigpf} 
from which  the notation is borrowed.
However, this second part is independent  of the first. Standard references for the language of derived categories
and constructible sheaves are  \ci{bo} and \ci{ks}; the reader may also consult \ci{sch}.

A variety is a separated scheme of finite type over the field of 
complex numbers. In particular, we do not assume that varieties are irreducible.
When dealing with intersection cohomology, it is convenient to 
assume irreducibility. However,  it is not necessary to do so; see
$\S$\ref{niv}.

The term stratification refers to algebraic Whitney stratifications. Algebraic
varieties and maps can be stratified.

Let $Y$ be a variety. 
We denote by ${\mathcal D}_Y$
the bounded constructible derived category, i.e.  the full subcategory of 
the derived category
$D (Sh_{Y})$
of the category of sheaves of abelian groups   on $Y$ consisting
of bounded  complexes $K$
 with constructible cohomology  
sheaves ${\mathcal H}^{l}(K)$ (though we do not pursue this, one can also get by with ``weak constructibility" \ci{ks}).     The objects are simply called constructible complexes.

The results of this paper that do not have to do with mixed Hodge theory, e.g.
the ones in $\S$\ref{tpfvso1} concerning the geometric description
of the perverse filtration,  hold if we replace sheaves of abelian groups
by sheaves of $R$-modules, where $R$ is a Noetherian commutative ring
with finite global dimension, e.g.  $\zed$, a principal ideal domain, a field, etc.

The variants  of these results  for varieties over an arbitrary field of definition 
and for the  various versions of  \'etale cohomology also hold, with very similar proofs,
and are left to the reader.

Given an algebraic map $f: X \to Y,$
we have the usual  derived functors $ f^{*},$ $f_{*}:= Rf_{*},$ $f_{!} : 
= Rf_{!}$
and $f^{!}$ acting between ${\mathcal D}_X$ and ${\mathcal D}_Y.$

Given  $K \in {\mathcal D}_Y$, we denote  the  (hyper)cohomology groups by  $H^*(Y, K)$
and the (hyper)cohomology groups with compact supports  
by $H_{c}^*(Y, K)$.

A $t$-category (cf. \ci{bbd, ks}) is  a triangulated category ${\mathcal D}$ endowed with a $t$-structure.
The
 truncation functors  are denoted   $\ts{a}:=\tu{a},$ $\tb{b}:=\td{b},$
the cohomology functors  $H\!\!:=$$ H^{0} \!\!:=\!\!\ts{0} \circ \tb{0},$ 
$H^{l}:= H^{0} \circ [l] = [l] \circ \ts{l}\circ \tb{l} =  [l] \circ  \tb{l}\circ  \ts{l}$
 have values in the abelian heart ${\mathcal C}$
of the $t$-category ${\mathcal D}$.

In this paper,  we deal with the standard and with the  middle-perversity 
$t$-structures on ${\mathcal D}_{Y}$. 
One word of caution: in the context of integer coefficients,
Verdier Duality does not preserve
middle perversity and there is no simple-minded exchange of
cohomology with cohomology with compact supports. Because of this,
at times we need to prove  facts in cohomology and in cohomology
with compact supports separately, though in our case the arguments 
run parallel, i.e. by inversion of the arrows. If one uses  field coefficients, these repetitions
can be avoided by invoking Verdier duality.

For the standard $t$-structure on ${\mathcal D}_Y$, we have the standard truncation functors,
the cohomology functors are the usual  cohomology sheaf functors ${\mathcal H}^*$ and the heart
is equivalent to the category of constructible sheaves
of $\zed$-modules on $Y$.
 
For the middle perversity
$t$-structure, the truncation and cohomology  functors are denoted by
$\ptu{l},$ $\ptd{l},$ $^{\frak p}\!{\mathcal H}^{l}$ and  the heart 
is the abelian category ${\mathcal P}_{Y} \subseteq {\mathcal D}_Y$
of (middle) perverse sheaves of $\zed$-modules on $Y$.

An  $n$-flag $Y_{*}$
(see Example \ref{sofclosed}) on a variety $Y$ is a
sequence of closed subvarieties of $Y:$
$$
\emptyset = Y_{-n-1} \subseteq Y_{-n} \subseteq Y_{-n+1} \subseteq \ldots
\subseteq Y_{-1} \subseteq Y_0 = Y.
$$ 
Typically, $Y_{-i}$ will be the intersection of $i$ hyperplane sections of 
a quasi projective $Y$ embedded in some projective space.
For technical reasons, the embedding must be chosen to be affine
(affine embeddings always exist). The  ``negative" indexing scheme
for the elements of a  flag serves the purposes
of this paper.

All filtrations on abelian groups  (and complexes) $M$ etc., are  decreasing, i.e. $F^{i} M \supseteq F^{i+1} M$ and  finite, i.e. $F^{\ll 0} M = M$
and $F^{\gg 0} M = 0$.
A  filtration 
is said to be of type $[a,b],$ if $F^a M =M$ and $F^{b+1}M=0.$

We  shall consider the following filtration  
 on cohomology as well as on cohomology with compact supports (see \S\ref{ssawts}):
the  standard $L_{\tau}$,  the Leray $L^f_{\tau}$,  the perverse 
$L_{\pp}$ and the  perverse Leray $L^f_{\pp}$ filtration.
The indexing scheme differs slightly from the one
in \ci{decmigpf}. This is to serve the purposes of this paper, especially the proof
of Propositions \ref{descrfiltr} and \ref{descrfiltr2}, where
the  indexing scheme employed in this paper places  conveniently the spectral sequences
in certain quadrants and  facilitates the  analysis.

Let $j: U \to Y \leftarrow Z : i$ be maps of varieties such that
$j$ is an open embedding and $i$ is the complementary closed embedding.
We have  distinguished triangles
and  long exact sequence of relative cohomology (in our situation we have that  $j^*=j^!$ and $i_!=i_*$)
\begin{equation}
\label{adt11}
j_!j^! K \to K \to i_* i^* K \stackrel{[1]}\to, \qquad
\ldots \to H^*(Y,K_U) \to H^*(Y,K) \stackrel{r}\to H^*(Z, K_{|Z}) \to \ldots
\end{equation}
\begin{equation}
\label{adt112}
i_!i^! K \to K \to j_* j^* K \stackrel{[1]}\to, \qquad
\ldots \to H^*_c(Z;i^!K) \stackrel{r'}\to H^*_c(Y,K) \to H^*_c(Y,Z, K) \to \ldots.
\end{equation}
The complex $K_U= j_! j^!K$ on $Y$ is not to be confused with $K_{|U}$ on $U$. Note also
that
$H^*(Y, K_U) = H^*(Y, Z, K)$.

The maps $r$ and $r'$ are called restriction and co-restriction maps.

In what follows, while the symbols are ambiguous, 
the formul\ae$\,$ are not. The slight abuse of notation is compensated
by the simpler-looking formul\ae.
Given a Cartesian diagram of maps of varieties
\begin{equation}
\label{cde}
\xymatrix{
X' \ar[r]^{g} \ar[d]^{f} & X \ar[d]^f \\
Y' \ar[r]^g & Y,
}
\end{equation}
there are
the base change maps
\begin{equation}
\label{bc1}
g^* f_* \lorw  f_* {g}^*, \qquad g^! f_! \longleftarrow f_! g^!
\end{equation}
and the base change isomorphisms
\begin{equation}
\label{bc2}
g^* f_! \simeq  f_! {g}^*, \qquad  f_* g^! \simeq g^! f_*.
\end{equation}

The base change maps (\ref{bc1}) are  isomorphisms if either
one of the following conditions is met:
\begin{itemize}
\item
 $f$ is proper; i.e. we have the
base change theorem for the proper
map $f$; 

\item $f$ is locally topologically trivial over $Y$;

\item $g$ is smooth; i.e. we have the base change theorem
for the smooth map $g$.
\end{itemize}

In fact: if $f$ is proper, then $f_!=f_*$ and (\ref{bc2}) implies (\ref{bc1}); similarly, if
$g$ is smooth of relative dimension $d$, since then  $g^! = g^*[2d]$ (cf. \ci{ks}, Proposition 3.2.3);
the remaining case follows easily from the K\"unneth formula.

The term ``mixed Hodge (sub)structure"   is abbreviated to MH(S)S.

\medskip
The following  is essentially due to M. Artin.
\begin{tm}
\label{cdav}
{\rm ({\bf Cohomological dimension of affine varieties})}
Let $Y$
be affine of dimension $n$,  $P \in \pe_{Y}$.
Then 
$$
H^{r}(Y, P) =0, \quad   \forall r \notin  [-n,0] 
\qquad  \qquad H^{r}_{c}(Y, P) =0, \quad  \forall r \notin [0,n].
$$
\end{tm}
{\em Proof.}  See  \ci{ks}, Proposition 10.3.3 and Theorem 10.3.8. 
See also
\ci{sch},  corollaries 6.0.3 and 6.0.4.
For the \'etale case, see \ci{bbd}, Th\'eor\`eme 4.1.1. 
\blacksquare

\begin{rmk}\label{usofrutto}
{\rm
In the context of this paper, a general hyperplane $H$ is one chosen as follows.
Pick any embedding of $Y$ into projective space ${\Bbb P}$; take the closure $\overline{Y}\subseteq {\Bbb P}$;
choose an algebraic  Whitney stratification $\Sigma$ of $\overline{Y}$ so that $Y$ is a union of strata;
choose, using the Bertini theorem, a
  hyperplane  $\overline{H}\subseteq {\Bbb P}$ so that it meets  transversally all the strata
  of $\Sigma$.  
Take $H: = Y \cap \overline{H}$. For a discussion see \ci{decbday}. See also \ci{decmigpf}, \S5.2.
}
\end{rmk}
\begin{tm}
\label{swl}
{\rm ({\bf   Lefschetz hyperplane theorem})}
Let $Y$ be quasi projective  of dimension $n$,   $P \in \pe_{Y} $.
Let  $H\subseteq Y $ be a general hyperplane section 
with respect to any embedding of $Y$ into projective space
and  $J: (Y \setminus H) \to Y$ be the corresponding open immersion.
We have:
$$
 H^{r}(Y, J_{!}J^{!}P) =0, \quad \forall r < 0,
\qquad \qquad
H^{r}_{c}(Y, J_{*}J^{*}P) =0, \quad \forall r > 0 .
$$
\end{tm}
{\em Proof.} This is due to several people:  Goresky and MacPherson 
\ci{gomasmt}, Deligne (unpublished) and Beilinson, \ci{be}, Lemma 3.3.
Deligne's and Beilinson's   proofs are  valid also in the \'etale case.   See also
\ci{sch}, p.397-398.

\blacksquare

\begin{rmk}
\label{wswl}
{\rm 
The  Lefschetz hyperplane theorem \ref{swl}  implies
that
\begin{enumerate}
\item the restriction map $H^r(Y, P) \to H^r(H, P_{|H})$ is an isomorphism
for every $r \leq -2$ and it is injective for $r= -1$; similarly for ${\mathcal F}$,
and
\item
the 
co-restriction map $H^r_c(H, i^!P ) \to H^r_c (Y, P)$ is an isomorphism
for $r \geq 2$ and surjective for $r=1$. 
\end{enumerate}
To my knowledge, there is no analogue of  2.  for  constructible sheaves ${\mathcal F}$.
Note also  that since $H$ is general, $i^*P[-1]  = i^!  P[1]$ (cf. \ci{decmightam}, Lemma 3.5.4.(b), 
or \ci{sch}, p.321)  is perverse on $H$.

\n
In \ci{decbday},  section 3.3.2, I incorrectly stated  that the Lefschetz hyperplane theorem
requires to first choose an affine embedding of $Y$ into projective space.

}\end{rmk}

\begin{rmk}\label{newrmk}
{\rm
For completeness, let us mention that  theorems \ref{cdav} and \ref{swl} admit well-known versions for 
constructible sheaves. Seeee \ci{decbday},  Appendix and  also \ci{sch}.
}
\end{rmk}

\medskip 
The following is due to  J.P. Jouanolou. 

\begin{pr}
\label{jo}
Let $Y$ be a quasi projective variety. There is a natural number $d$ and  a Zariski locally trivial
${\Bbb A}^{d}$-fibration $\pi: {\mathcal Y}\to Y$ 
with affine transition functions and affine total space  ${\mathcal Y}$.
\end{pr}
{\em Proof.} See \ci{jou}. See also \ci{decbday}.  
\blacksquare

\begin{rmk}
\label{d=n}
{\rm
There is no canonical choice for the fibration. One can arrange for $d= \dim{Y}$
but, in general, not less, e.g. $Y= {\Bbb P}^n.$
}
\end{rmk}

\section{The perverse filtrations via spectral objects}
\label{tpfvso1}
The goal of this section is to prove the results in $\S$\ref{tgofsss}. This is done
using the preparatory results on spectral objects and spectral sequences in $\S$\ref{sezodi}.
Spectral objects are discussed in the next $\S$\ref{subspob}.

\subsection{Spectral objects and spectral sequences}
\label{subspob}

\subsubsection{Spectral objects in a triangulated category}
Spectral objects have been introduced by Verdier. A good reference is
\ci{deseattle}.

Let ${\mathcal D}$ be a triangulated category.
A {\em spectral object}  in  ${\mathcal D}$ is the data:

\smallskip
\n
$(1)$ a family of objects $X_{pq}$ of ${\mathcal D}$ indexed by pairs $p\leq q \in \zed$,

\smallskip
\n
$(2)$ for $p'\leq p,$ $q' \leq q,$ a morphism $X_{pq} \to X_{p'q'}$,

\smallskip
\n
$(3)$  for $p\leq q \leq r,$ a morphism $\partial:X_{pq} \to X_{qr}[1]$ 
called {\em coboundary}

\smallskip
subject to the  requirement that

\smallskip
\n
$(a)$ the morphisms (2) define a {\em contravariant} functor from the
category of ordered pairs
$(p,q),$ with $p \leq q,$ to ${\mathcal D}$,

\smallskip
\n
$(b)$ for $p\leq q\leq r,$ $p'\leq q' \leq r',$ and $p'\leq p,$  $q'\leq q,$ $r' \leq r,$
 the diagram
$$
\xymatrix{
    X_{pq} \ar[r]^{\partial} \ar[d] & X_{qr}[1]\ar[d] \\
      X_{p'q'} \ar[r]^{\partial} & X_{q'r'}[1]}
$$
of morphisms from (2) and (3) is commutative, 

\smallskip
\n
$(c)$ for $p\leq q\leq r$ the triangle
$$
X_{qr} \lorw X_{pr} \lorw X_{pq} \stackrel{\partial}\lorw X_{qr}[1]
$$
is distinguished.

\smallskip
There is an obvious notion of morphism of spectral objects and  spectral objects
in ${\mathcal D}$ form a category.

The axioms imply $X_{pp}=0$ so that, on the $(p,q)$-plane,  the display
occurs on and above the line $q = p+1.$

We  work exclusively  with {\em bounded} spectral objects, i.e.
objects for which
$X_{p,p+1} =0$ for  $|p| \gg 0$. If $X_{p,p+1}=0$ for
  $p \neq [a,b]$, then we  
we say the spectral object has   {\em amplitude in } $[a,b].$
In this case: 
i)  all $X_{pq}$ with $p <a$ and $q>b$ are isomorphic  to one another
via the maps
(2)
and we denote them  by
$X_{-\infty, \infty},$  ii)  all the $X_{p,q},$ with $p$ fixed and $q >b,$ are 
 isomorphic via the maps (2)
and we denote them by $X_{p, \infty}.$

Up to applications of the translation functor $C \mapsto C[1]$, one may choose to
consider only the case of amplitude in
$[0,b-a],$ or $[a-b, 0].$
  In the former case we have that $X_{0,q} \simeq
X_{-1,q} \simeq X_{-2,q} \ldots,$ for every $q$  and the essential part 
of the display is a triangle in the first quadrant.

Giving a spectral object with
amplitude in  $[0,1]$ is the same as giving  a distinguished triangle. 
Amplitude in  $[0,2]$ corresponds to an octahedron diagram 
(as in the octahedron axiom for triangulated categories), etc.
A spectral object is  therefore a suitably compatible system of triangles, octahedrons, etc.

\begin{ex}
\label{exfc}
{\rm ({\bf Filtered complexes})
Let ${\mathcal A}$ be an abelian category and 
$(K,F)$  be a filtered complex with finite filtration. 
By setting  $X_{pq}:= F^pK/F^q K,$ one gets 
 a bounded  spectral object
in the derived category  $ D({\mathcal A})$.  In other words, 
an object of the finite filtered derived category $DF({\mathcal A})$
yields a spectral object in $D({\mathcal A})$.
}
\end{ex}

\begin{ex}
\label{somor}
{\rm ({\bf Sequence of maps})
Let  $\mathcal A$ be an abelian category with enough injectives and
\[
 \ldots \lorw K_{i+1} \lorw K_{i}  \lorw  K_{i-1} \lorw \ldots \]
be a sequence of morphisms in $D^+({\mathcal A}),$ 
 with $K_i =0,$ for $i \gg 0$
and
 $K_{i+1} \simeq K_i,$
for $i \ll 0.$  
There exists a filtered complex, with finite filtration, 
$(K,F)$ in $D^+({\mathcal A}),$ such that the given sequence
is isomorphic to the sequence of sub-complexes $F^i K$
(see \ci{bbd}, p.77).  In the cases we use in this paper,
this correspondence is canonical,
i.e. defined up to unique isomorphism in the filtered derived category.
In view of Example \ref{exfc}, we shall
speak of the spectral object in $D^+({\mathcal A})$  associated with a
sequence of morphisms as above.
}
\end{ex}

\begin{ex}
\label{socts}
{\rm ({\bf Truncation})
Let ${\mathcal D}$ be a $t$-category with truncation functors
$\tb{i}, \ts{i}$ and cohomology functors $H^i = [i] \circ \ts{i} \circ \tb{i}$.
An object $X$ of ${\mathcal D}$ yields a spectral object in ${\mathcal D}$ by setting:
$$
X_{pq}: = \ts{-q+1}\tb{-p}X.
$$
This correspondence is functorial.
If $X$ is bounded, then $X_{p,\infty} = \tb{-p}X.$
This choice  of indexing leads to  a decreasing filtration in cohomology (cf. 
$\S$\ref{assss}.(\ref{filtf})).
One has
$$
X_{p,p+1} =   \ts{-p}\tb{-p}X = H^{-p}(X)[p].
$$
Cohomological amplitude
$[0,n]$ here means that $H^{l}(X) =0$ for $l \notin  [-n,0].$ 
}
\end{ex}

\begin{rmk}
\label{spcssocts}
{\rm
The category ${\mathcal D}_Y$  of  bounded constructible complexes on a variety $Y$
admits the  standard  $t$-structure, i.e.
the one associated with the natural truncation functors. It also admits
the  distinct $t$-structures  associated with
distinct perversities.  
Each $t$-structure yields  spectral objects as in  Example \ref{socts} which
are isomorphic  
to the spectral objects we obtain by virtue of Example \ref{somor}
applied to the sequence of  truncation maps associated with the $t$-structure:
$\;
\ldots \tb{- i-1 } X \lorw \tb{ -i} X \lorw  \ldots  \;
$
}
\end{rmk}

\begin{ex}
\label{sofclosed}
{\rm {\bf (Sequence of closed subvarieties})
A {\em $n$-flag $Y_{*}$ } on the variety $Y$ is defined to be a sequence 
$$
\emptyset = Y_{-n-1} \subseteq Y_{-n} \subseteq Y_{-n+1} \subseteq \ldots
\subseteq Y_{-1} \subseteq Y_0 = Y
$$
of   closed subvarieties of $Y$. The indexing scheme  is chosen
as to serve the needs of this paper and 
obtain filtrations $F_{Y_{*}}$
of type $[-n,0]$ in cohomology and of type $[0,n]$
in cohomology with compact supports. 
Let
$$
j_p : Y \setminus Y_{p}  \lorw Y, \qquad
i_p: Y_p \lorw Y, \qquad
k_p: Y_p\setminus Y_{p-1} \lorw Y 
$$
be the corresponding  embeddings.

 Let $K \in {\mathcal D}_Y$.   Note that 
 ${j_p}_! j^*_p K = K_{Y-Y_{p-1} }$ and that
 ${i_p}_* i_p^! K = R \Gamma_{Y_p} K$.

 Without loss of generality, we assume that
 $K$ is injective so that 
 $R\Gamma_{Y_p} K = \Gamma_{Y_p} K$; this fact allows us, in particular, to form the quotients
 in (\ref{33}).

We have the  two spectral objects $K_{pq}$ and $K^!_{pq}$
 associated with the two  filtered complexes
\begin{equation}\label{somsof}
   (K, F_{Y_*}), \qquad 
F^p_{Y_*} K:= K_{Y-Y_{p-1}}, \qquad \qquad
K_{pq}: = F^p_{Y_*}  K /F^q_{Y_*}K,
\end{equation}
\begin{equation}
\label{33} 
  (K, G_{Y_*}), 
 \qquad 
G^p_{Y_*} K: = R\Gamma_{Y_{-p}} K, \qquad \qquad
K^!_{pq}: = G^p_{Y_*} K /  G^q_{Y_*} K.
\end{equation}
The spectral object $K_{pq}$
has amplitude in the interval $[-n,0]$ and $K^{!}_{pq}$ in $[0,n]$
and the resulting filtrations have type $[-n,0]$ and $[0,n]$, respectively.
}
\end{ex}

\subsubsection{Spectral objects in an abelian category with translation functor}
Let $\mathcal B$ be an abelian category and $\mathcal A$ be the associated category of graded objects.
The abelian category $\mathcal A$ is naturally endowed with the translation functor,
denoted $[1],$ which is an autoequivalence.

A spectral object in  $\mathcal A$ is defined analogously to one in ${\mathcal D}$ with the following
adapted axioms:

$(1')$ the objects $X_{pq}$ are in fact collections $\{X_{pq}^n\}_{n \in \zed}$ 
in ${\mathcal B};$

$(3')$ $\partial: X_{pq}^n \to X_{qr}^{n+1};$

$(c')$ there  is a long exact sequence
$$
\ldots \lorw X_{pr}^n \lorw  X_{pq}^n \stackrel{\partial}\lorw X_{qr}^{n+1} 
\lorw X_{pr}^{n+1} \lorw \ldots
$$

Let ${\mathcal D}$ be a  triangulated category, $\mathcal B$ be an  abelian
category and 
$T:{\mathcal D} \to {\mathcal B}$ be a cohomological functor.
Denote by $T^{\bullet}:{\mathcal D} \to {\mathcal A}$ the induced functor:
$T^n(X):= T(X[n]).$

In either of the categories ${\mathcal D}$ and $\mathcal A,$
there is the obvious notion of morphism of spectral objects 
$X_{pq} \to X'_{pq}.$
Spectral objects in ${\mathcal D}$ and in ${\mathcal A}$ form categories and 
$T^{\bullet}$ is a functor transforming spectral objects in ${\mathcal D}$ 
into spectral objects in ${\mathcal A}:$ if $X$ is a spectral object in ${\mathcal D},$ then
it is immediate to verify that
we obtain a spectral object $T(X)$
  in $\mathcal A$ by setting
$$
  T(X)^n_{pq} := T^{n}(X_{pq}).
$$

\subsubsection{Spectral objects and  spectral  sequences}
\label{assss}
Let $T : {\mathcal D} \to {\mathcal B}$ and $T^{\bullet}: {\mathcal D} \to {\mathcal A}$  be as above. 
Then $T$ 
also induces a functor
into the category of spectral sequences. 
If $X $ is a bounded spectral object in ${\mathcal D},$ then there is the object 
$X_{-\infty, \infty},$ canonically isomorphic to $X,$ and the
spectral sequence $E_1 (T(X)):$
\begin{equation}
\label{aspp}
 \qquad E_1^{pq} (T(X)) :=   T^{p+q} (X_{p,p+1})  \Longrightarrow 
T^{p+q} (X)
:= T^{p+q} (X_{-\infty, \infty} ).
\end{equation}
The associated  decreasing and finite filtration is 
\begin{equation}
\label{filtf}
F^p T^{u}(X) = \im \, \{ T^{u}(X_{p,\infty}) \}
\subseteq T^{u}(X), \qquad
F^p T^{p+q} / F^{p+1} T^{p+q}  = E_{\infty}^{p, q}.  
\end{equation}

\subsubsection{Spectral sequences:
re-numeration,  the $L$ filtration, translation}
\label{ssecss}
Let
$
E_{1}^{pq} \Longrightarrow
T^{p+q}
$
be a spectral sequence as in $\S$\ref{assss}
with abutment the filtration $F$. 

In view of Propositions \ref{descrfiltr} and \ref{descrfiltr2},
we introduce a different
indexing scheme.
The {\em re-numeration} ${\mathcal E}_{r+1}$  of $E_r,$ 
is the {\em same} spectral sequence,  with a different
indexing scheme:
\begin{equation}
\label{ren}
{\mathcal E}^{st}_{r+1}  := E_r^{-t, s+2t}, \quad
r\geq 1, \qquad 
{\mathcal E}_2^{st} \Longrightarrow T^{s+t}.
\end{equation}
One still has the  filtration $F,$ but also has the 
{\em $ L$-filtration} (which we adopt in this paper)
\begin{equation}
\label{lf}
F^p T^u =: {L}^{p+u}T^u, \qquad
L^{s}T^{s+t}/L^{s+1}T^{s+t} = {\mathcal E}^{st}_{\infty}.
\end{equation}

The following relations are  easily verified: 
\begin{equation}
\label{proppfl}
F^p ( T^{u} (X[d]) ) = F^{p-d} ( T^{u +d } (X) ),
\qquad
L^p (T^{u} (X[d]) ) = L^p ( T^{ u  +d } (X) ).
\end{equation}

\begin{ex}
\label{grren}
{\rm If one re-numbers the usual Grothendieck spectral sequence
for cohomology
$ E_{1}^{pq}= H^{p+q}(Y, {\mathcal H}^{-p}(K)[p]) \Longrightarrow H^{p+q}
(Y, K),$
which has 
$${F}^{p}H^j(Y,K)= \im\, \{ H^*(Y,\tau_{\leq -p}K )  \lorw H^*(Y,K)\},$$ 
then 
one obtains
${\mathcal E}_{2}^{st}= H^{s} (Y,{\mathcal H}^{t} (K))
\Longrightarrow H^{p+q}
(Y, K)$
 with decreasing $L$-filtration
$$L^{s}H^{j} (Y,K)= \im\, \{  H^{j}(Y,\tau_{\leq -s+j} (K)\}.$$
Note that the re-numbered Grothendieck spectral sequence
is in the quadrants I and IV.
}
\end{ex}

The {\em $l$-translate} $G(l)$ of a filtration $G$
is defined by setting
 \begin{equation}
 \label{deftf} 
 G(l)^s:= G^{l+s}.
\end{equation}
The {\em $l$-translate} $E_r(l)$ of a  spectral sequence
$E_r$ with abutment a filtration $F$ is defined as follows:
\begin{equation}
\label{tsspe}
E_{r}(l)^{pq}: = E_{r}^{p+l, q-l} \Longrightarrow T^{p+q}.
\end{equation}
The  associated  filtrations satisfy,
\begin{equation}
\label{trfi}
F (E_{1}(l)) = F(E_{1}) (l).
\end{equation}
This formula  
holds also for the associated  $L$ filtration and one has
$F(l)^{p} T^{u} = L(l)^{p+u}T^{u}.$

\subsubsection{The spectral sequence associated with a $t$-structure:
standard and perverse spectral sequences and filtrations} 
\label{ssawts}
 Let ${\mathcal D}$ be a  $t$-category, let
$T: {\mathcal D} \to {\mathcal B}$ be a cohomological functor and  let  $T^{\bullet}:{\mathcal D} \to 
{\mathcal A}$ be the
associated graded functor. 
Given $X \in {\mathcal D},$
denote by $X$ also the associated spectral object  via Example \ref{socts}.
Note that $T(X)$ has now two different meanings: $T(X)  \in {\mathcal B}$ 
and $T(X)$ the  spectral object
in ${\mathcal A}.$
I hope this does not generate confusion.

Recall the formula $X_{p,p+1}= H^{-p}(X)[p]$.
Assume that
 $H^l(X)=0,$ for $|l|\gg 0,$ i.e. that  the spectral object $X$ 
is bounded. There is the spectral sequence $E_1 (T(X)) $ as in
$\S$\ref{assss}:
$$
 E_1^{pq}:=   
T^{p+q} (H^{-p}(X)[p]) = T^{2p+q} (H^{-p}(X)) \Longrightarrow T^{p+q} (X).
$$
If we re-number as in  (\ref{ren}),  then we have
${\mathcal E}_2 (T(X))$:
$$
 {\mathcal E}_2^{s,t} T^s(H^t(X)) 
 \Longrightarrow  T^{s+t}(X).
$$
For the associated filtrations, we have
$
F^p T^{u}(X) =
\im \, \{  T^{u} (\tb{-p} X) \to T^{u}(X)\} 
\subseteq T^{u}(X), 
$ and we have  the equalities (\ref{lf}).

In this paper,  we work with    the standard and with  the middle perversity
$t$-structures on ${\mathcal D}_Y$, denoted by
 $\tau$ and $\pp$ respectively,
and with the cohomology and compactly supported cohomology
functors, i.e.  with the special cases when $T= H^0(Y,-)$, or $T=H_c^0(Y,-)$. 

Let $Y$ be a  variety 
and $K \in {\mathcal D}_Y.$ We have the following spectral sequences
with associated filtrations:
\begin{equation}
\label{gro}
{\mathcal E}_2^{st} (K, \tau):= H^s(Y, {\mathcal H}^t (K)) \Longrightarrow
(H^{s+t}(Y, K), L_{\tau}),
\end{equation}
\begin{equation}
\label{per}
{\mathcal E}_2^{st} (K, \pp):= H^s(Y, \pc{t}{K}) \Longrightarrow
(H^{s+t}(Y, K), L_{\pp}).
\end{equation}
We have analogous sequences 
$_c{\mathcal E}$ for compactly supported cohomology $H^*_c(Y,K).$

The  spectral sequence (\ref{gro}) is called the {\em standard} (or Grothendieck) {\em spectral sequence}
and the associated filtration  of type $L$ (versus type $F$)   is called the
 {\em standard filtration} and is denoted $L_\tau$.
The spectral  sequence (\ref{per}) is called the {\em perverse spectral sequence}
and the associated filtration of type $L$
is called the  {\em  perverse filtration}
 and is denoted by $L_{\pp}$.

Let $f: X \to Y$ be a map of varieties and let $C \in {\mathcal D}_X.$ 
We have the Leray and the perverse Leray spectral sequences
with associated Leray and perverse-Leray filtrations:
\begin{equation}
\label{ler}
{\mathcal E}_2^{st}(f_{*} C, \tau) \Longrightarrow
(H^{s+t}(X, C), L_{\tau}^f),  \qquad
_c{\mathcal E}_2^{st}(f_{!} C, \tau) \Longrightarrow
(H^{s+t}_c(X, C), L_{\tau}^f);  
\end{equation}
\begin{equation}
\label{pler}
{\mathcal E}_2^{st}(f_{*} C, \pp) \Longrightarrow
(H^{s+t}(X, C), L_{\pp}^f),  \qquad
_c{\mathcal E}_2^{st}(f_{!} C, \pp) \Longrightarrow
(H^{s+t}_c(X, C), L_{\pp}^f).  
\end{equation}
\begin{rmk}
\label{tipi}
{\rm
The standard spectral sequences are in the quadrants
I and IV.
Let $Y$ be affine
of dimension $n$. Theorem \ref{cdav} implies that
the standard filtration is of type
$[0,n]$  in cohomology and of type  $[n,2n]$ in cohomology with compact supports
and the corresponding standard spectral sequences
${\mathcal E}_{r}^{s,t}$, $_c{\mathcal E}_{r}^{s,t}$ live in the columns $s \in [0,n]$ and $s \in [n,2n]$, respectively.
The perverse filtration has types $[-n,0]$ and $[0,n],$ respectively,
and  the  perverse  spectral sequences are in the quadrants
II and III (cohomology), with columns $s\in [-n,0]$, and I and IV (cohomology with compact supports), with columns $s\in [0,n]$.
If $Y$ is a variety of dimension $n$ covered by
$1+t$ affine open sets, then the four types  for the filtrations above are
$[0, n+t], $ $[n-t, 2n],$ $[-n,t]$ and $[-t, n],$ respectively and the spectral sequences
live in the corresponding columns.}
\end{rmk}

\subsubsection{The spectral sequence associated with a flag  of closed subspaces.}
\label{ssflag}
Let $Y$ be a variety, let $K \in {\mathcal D}_Y$ 
and  let $Y_{*}= \{ \emptyset \subseteq Y_{-n} \subseteq 
 \ldots \subseteq Y_{0} =Y \}$ 
be a $n$-flag
as in Example \ref{sofclosed}. There are the  two spectral objects 
$K_{p,q}$ and $K^{!}_{pq}$
in ${\mathcal D}_Y$. 
Their amplitudes are  in $[-n,0]$ and in $[0,n], $ respectively.
If we apply the cohomology functors $H^*(Y,-)$ to the first one and
$H^*_c(Y,-)$ to the second one, then   we get two spectral objects
whose associated spectral sequences are the classical 
spectral sequences for the  filtration of $Y$ into closed subspaces:
\begin{itemize}
\item
the spectral sequence for cohomology and the associated filtration are \,
$E_{1} (K, Y_{*}):$
\begin{equation}
\label{coss}
E_1^{p,q} = H^{p+q}(Y, K_{Y_p- Y_{p-1} } ) \Longrightarrow
H^{p+q}(Y,K),
\end{equation}
\begin{equation}
\label{filtfl}
F^p_{Y_{*}} H (Y, K) = 
\ker \, \{ H^*(Y, K) \lorw H^*(Y_{p-1}, K_{p-1} )     \};
\end{equation}
\item
the spectral sequence and filtration for cohomology with compact supports
are
$_{c} E_{1}(K, Y_{*}):$
\begin{equation}
\label{cosan}
_cE_1^{p,q} = H^{p+q}_c(Y,  R\Gamma_{Y_{-p} -Y_{-p-1}}K) \Longrightarrow
H^{p+q}_c(Y,K),
\end{equation}
\begin{equation}
\label{444}
G^p_{Y_{*}} H_c^* (Y, K) = \im \, \{  H^*_c(Y_{-p}, i_{-p}^!K )
\lorw H^*_c(Y,K) \}.
\end{equation}
\end{itemize}

\begin{rmk}
\label{tipfiltr}
{\rm
Note that $F_{Y_*}^{-n} H^*(Y, K)= H^*(Y,K)$ and $F^1 H^*(Y,K) =0,$
i.e. the filtration  $F_{Y_{*}}$ in cohomology has type $[-n,0].$
Analogously,  $G_{Y_*}$  has type $[0,n]$.
}
\end{rmk}

\subsection{Preparatory results}
\label{sezodi}

\subsubsection{Translating spectral objects}
\label{trspob}
A map $X\to Y$  in a   $t$-category   ${\mathcal D}$ yields
 a map of spectral objects $X_{pq} \to Y_{pq}.$ By applying 
a cohomological functor $T:{\mathcal D} \to {\mathcal B}$ and by recalling
the re-numeration
({\ref{ren}),
we get a morphism of spectral sequences 
$$
{\mathcal E}_r^{st}(X) \lorw {\mathcal E}_r^{st}(Y), \quad r \geq 2.
$$
\begin{ex}
\label{triv}
{\rm 
Let $i: pt \to {\Bbb A}^1.$ We have the adjunction morphism
$$
\rat_{ {\Bbb A}^1}[1] \lorw i_* i^* \rat_{ {\Bbb A}^1}[1]  =\rat_{pt}[1] .
$$
We consider the standard and the middle-perversity $t$-structures.
The map of spectral sequences associated with the standard $t$-structure
is an isomorphism of spectral sequences and it induces an isomorphism
on the  filtered abutments.
On the other hand,  middle-perversity  yields the zero map of spectral sequences.
Moreover,  we have
\[ F^0 H^{-1}({\Bbb A}^{1},\rat_{ {\Bbb A}^{1}}   [1]) = \rat, \quad 
F^1 H^{-1}({\Bbb A}^{1}, \rat_{ {\Bbb A}^{1}}   [1]) =0, \quad Gr_F^0 = \rat, \quad Gr_F^1 = 0,
\]
\[ F^1 H^{-1}(pt,\rat_{ {pt}}   [1]) = \rat, \quad 
F^2 H^{-1}(pt , \rat_{ {\Bbb A}^{1}}   [1]) =\rat, \quad Gr_F^0 = 0 \quad Gr_F^1 = \rat,
\]
Similarly, for the $L$ filtrations for which we have, respectively:
$Gr_L^{-1} = \rat$, $Gr_L^0 = 0$  and
$Gr_L^{0} = \rat$, $Gr_L^1 = 0$.  In particular,  a mere re-numbering
of the spectral sequences will not yield an isomorphism.
The reason for this different behavior is explained by 
Lemma \ref{sotex}: the functor  $i^{*}$ is $t$-exact for the
standard $t$-structure, but it is not for the middle-perversity $t$-structure.
}
\end{ex}

In the  example above   one could re-index the filtration on $H^{-1}(pt, \rat_{pt} [1])$ so
that, for example,
$
L^{-1} H^{-1} (\rat_{{\Bbb A}^1} [1]) \simeq
L^{-1} H^{-1}(\rat_{pt}[1]).
$
One can do so, but may also care to define a 
morphism of spectral sequences compatible with this objective. 
In this case, this is possible, due to the presence of suitable 
$t$-exactness
for $i^{*}[-1]$ with respect to $\rat_{ {\Bbb A}^1} [1];$
see Remark \ref{1122}.

The necessity of considering a translated filtration on, say, the second spectral sequence
while at the same time retaining a map of spectral sequences,
 has lead  us to
the notion  of translated spectral objects.  We  need to have a meaningful
map of spectral sequences  after this translation has taken place. Lemma
\ref{sotex} is a criterion for reaching this objective.

Recall the notions  (\ref{deftf}) of translated filtration,
 (\ref{tsspe}) of  translated spectral sequence  and the effect
 (\ref{trfi}) of translation
  on the abutted  $F$ and $L$ filtrations.

\begin{defi}
\label{deftrso}
{\rm 
Let $X = \{ X_{pq} \}$ be a spectral object in ${\mathcal D}$  ($X=\{ X^n_{pq} \}$ 
be a spectral object in $ {\mathcal A},$ respectively). Let $l \in \zed$ and
define
the {\em translation of $X$ by $l$} to be  the spectral object
$X(l)$:
$$
X(l)_{pq} := X_{p+l, q+l}, \qquad
(X(l)^n_{pq} := X^n_{p+l, q+l},\quad \mbox{respectively}).
$$
}
\end{defi}

In the presence of a cohomological  functor $T:{\mathcal D} \to {\mathcal B},$ a spectral object
$X$ in ${\mathcal D}$ gives the spectral object $T (X)$ in $\mathcal A$ and we have
$
(T(X)(l))^n_{pq} = T^{n}(X_{p+l, q+l}).
$


We have the following simple equalities
\begin{equation}
\label{tuttook}
E_{1} (  T(X(l))  ) = E_{1} (T(X)) (l) \Longrightarrow
(T(X), L(l)).
 \end{equation}

\bigskip
Let   $d \in \zed$  and $v:{\mathcal D} \lorw {\mathcal D}'$ be a
functor of $t$-categories, i.e. $v$ is additive, it commutes with translations
and it sends distinguished triangles into distinguished
triangles.  We say that  $v\circ [d] $ is $t$-exact if   
\begin{equation}
\label{serve}
v ( \ts{a} (X)) = \ts{a+d} ( v(X) ), \qquad
v ( \tb{a} (X)) = \tb{a+d} ( v(X) )
\end{equation}
or, equivalently, if:
$$
[ v \circ  [d], \ts{a} ] = [ v\circ [d], \tb{a} ]  =[ v \circ [d], H^a ]=0.
$$

\begin{lm}
\label{sotex}
Let
$$
\xymatrix{
{\mathcal D} \ar[rr]^{u^{\bullet}}  \ar[dr]_{T} & & {\mathcal D}'  \ar[dl]^{T'} \\
& {\mathcal B} &
}
$$
be a  commutative diagram of functors, with ${\mathcal D}$ and ${\mathcal D}'$ $t$-categories
and ${\mathcal B}$ abelian
such that: 

\begin{enumerate}
\item
$u^{\bullet} [d]$ is $t$-exact, for a fixed $d \in \zed,$

\item$T$ and $T'$ are  cohomological,

\item
there exists $u_{\bullet}: {\mathcal D}' \to {\mathcal D}$
such that  $(u^{\bullet}, u_{\bullet})$ 
is an adjoint pair  and

\item
 $T \circ u_{\bullet} = T';$
in particular,   $T u_{\bullet} u^{\bullet} = T' u^{\bullet}.$
\end{enumerate}
 Let $X \in {\mathcal D}$ and  $X_{pq}= \ts{-q+1} \tb{-p}X$ be the associated spectral object
as in Example {\rm \ref{socts}}.
Then there is a morphism of spectral objects in ${\mathcal A}$:
$$
T(X) \lorw  [{T'} ( u^{\bullet} X)](-d)
$$
yielding a morphism  of spectral sequences
$$
E_{1} (T(X))  \lorw  E_{1} (T' (u^{\bullet} X)(-d))
$$
with resulting filtered map
$$
(T^{p+q}(X),L) \lorw  ({T'}^{p+q}(u^{\bullet}X), L(-d)).
$$
Analogously, if  $(u_{\bullet}, u^{\bullet})$ is an adjoint pair,
then the statement with arrows reversed  holds true.
\end{lm}
{\em Proof.} We prove the first statement. The proof of the second one is identical.

By adjointness, there is the morphism of functors 
$a: Id  \to u_{\bullet} u^{\bullet}$
inducing a map
$$
T^n( \ts{-q+1} \tb{-p} X ) \lorw 
T^n (u_{\bullet} u^{\bullet}\ts{-q+1} \tb{-p} X) \stackrel{4.}=
{T'}^n (u^{\bullet}   \ts{-q+1} \tb{-p} X ).
$$

On the other hand, by $t$-exactness
(\ref{serve}) and by the definition of translation \ref{deftrso},
we get the  natural identifications
$$
{T'}^n (u^{\bullet}   \ts{-q+1} \tb{-p} X ) =
{T'}^n (   \ts{-q+d+1 } \tb{-p+d}   u^{\bullet} X ) =
[   {T'}( u^{\bullet}  X) (-d)]^{n}_{pq}.
$$

We thus get a natural map induced by the adjunction map $a: X \to u_{\bullet}
u^{\bullet}X$:
$$
T(X)_{pq}^n =
T^n( \ts{-q+1} \tb{-p} X ) \lorw {T'}^n (   \ts{-q+d+1 } \tb{-p+d}   u^{\bullet} X ) =
T'(u^{\bullet} X)(-d)_{pq}.
$$

Since $a$ is a morphism of functors, checking that the map above induces
a morphism of spectral objects in ${\mathcal A}$ reduces to a formal verification.
\blacksquare

\begin{rmk}
\label{itsdiff}
{\rm
Example \ref{triv} shows that the map of spectral sequences
arising via functoriality from  the adjunction map $X \to u_{\bullet} u^{\bullet} X$
is {\em not} the one of Lemma
\ref{sotex} which needs the extra input of $t$-exactness. Note also
that the realigned map of spectral sequences $E_1 \to E_1 (-d)$ does not
correspond to a  realigned map of spectral objects $X \to u_{\bullet} u^{\bullet} X [-d]$.  
}
\end{rmk}

\begin{rmk}
\label{expl}
{\rm 
The map of spectral sequences  of Lemma \ref{sotex}
is, a bit more explicitly:
$$
T^{p+q} (H^{-p}(X)[p])  \lorw {T'}^{p+q}( H^{-p+d}(u^{\bullet}X)[p-d] ).
$$
}
\end{rmk}

\begin{rmk}
\label{1122}
{\rm
The   functor $i^{*} [-1]$ of Example \ref{triv} is  is not $t$-exact
for the middle perversity $t$-structures
since, for example, $i^{*} \rat_{pt} [-1]$ is not perverse. However, 
it is $t$-exact when applied to complexes, such as 
$\rat_{ {\Bbb A}^1} [1]$ which are constructible with respect to a stratification
for which the inclusion $i: pt \to  {\Bbb A}^1$
is   a normally nonsingular inclusion of codimension one. 
Such complexes form a $t$-subcategory and $i^{*} [-1]$ is 
$t$-exact when acting on such complexes.
It follows that one can apply Lemma \ref{sotex} with $u =i^{*},$ $d=\!\!-1$  
and get the desired isomorphism of spectral sequences.
}
\end{rmk}

\begin{rmk}
\label{flss}
{\rm
Lemma \ref{sotex} implies the usual functoriality of the Leray spectral sequence
with respect to square  commutative diagrams (\ref{cde}). One only needs to
observe that $g^*$ is exact in the usual sense, i.e. $t$-exact for the
standard $t$-structure. In fact, we have the adjunction map
$f_* C \to g_* g^* f_*C$ and the natural map $g^* f_* C \to f_* g^* C$.
By applying   the lemma to the first and   functoriality to the second, we obtain
$H^s(Y, R^t f_* C) \to H^s (Y', {\mathcal H}^t(g^* f_* C) ) \to H^s(Y', R^sf_*(g^*C))$.
}
\end{rmk}

\subsubsection{Inductive behavior of spectral sequences}
\label{alossq}
In this section we prove Proposition
\ref{descrfiltr} and its analogue  Proposition \ref{descrfiltr2}.
These propositions are used in the proof of Theorems
\ref{mtm} and \ref{flplr}, i.e. the proof of the   geometric description
of the perverse and perverse Leray filtrations  in the  quasi projective case.

The following  lemma is spelled-out for the reader's convenience
\begin{lm}
\label{pippa}
Let 
$$
\xymatrix{
    E'\ar[r]^{d'_E} \ar[d]^{\phi'}  &  E \ar[r]^{d''_E} \ar[d]^\phi  & E'' \ar[d]^{\phi''} \\
      F' \ar[r]^{d'_F} &   F \ar[r]^{d''_F} & F''
      }
$$
be a  commutative  diagram in an Abelian category such that
$d^2=0, $ $ \phi'$ epic, $\phi$ monic ($\phi'' $ monic, $\phi$ epic, respectively). 
Then the induced map on $d$-cohomology
$$
H_E \stackrel{H^*({\phi})}\lorw H_F
$$
is monic (epic, respectively).
\end{lm}
{\em Proof.} 
Since $\phi$ is monic, the induced map on kernels is monic and since $\phi'$ is epic,
the $d'$-images map isomorphically onto each other. This proves the first
statement.
As to the second one, 
since $\phi$ is epic and $\phi''$ is monic, the induced map on
 kernels is epic and we are done.
\blacksquare

\bigskip
The following Lemma is the key step in the proof
of the propositions that follow. The first statement is used to study cohomology
and the second one to study cohomology with compact supports.
We state it in the way we shall use it, i.e. using the
$({\mathcal E}, L)$-notation of (\ref{ren}) and (\ref{lf}),
 i.e.  $L^s H^{s+t}/L^{s+1}H^{s+t} = {\mathcal E}^{s,t}_{\infty}$

\begin{lm}
\label{ssis}
$\,$

\begin{enumerate}
\item
Let $\phi : {\mathcal E} \to {\mathcal F}$ be a morphism of bounded {\rm II/III}-quadrants spectral sequences with abutments $H_{\mathcal E}$ and $H_{\mathcal F}.$
Assume that $\phi_2^{pq} : {\mathcal E}_2^{pq} \to {\mathcal F}_2^{pq}$ is an isomorphism for every $p\leq -2,$ injective
for $p=\!\!-1$ and zero for $p=0.$ 
Then $\phi_\infty^{pq}$ is injective for every $p\leq -1$ 
and
$$\ke{\, \{H^*(\phi):  H_{\mathcal E} \to H_{\mathcal E}}\}  =  L^0 H_{\mathcal E} .$$

\item
Let  $\phi : {\mathcal E} \to {\mathcal F}$ be a morphism of bounded {\rm I/IV}-quadrants spectral sequences.
Assume that $\phi_2^{pq} : {\mathcal E}_2^{pq} \to {\mathcal F}_2^{pq}$ is an isomorphism for every $p>1,$ surjective
for $p=1$ and zero for $p=0.$
Then $\phi_\infty^{pq}$ is surjective for every $p\geq 1$
and
 $$L^1 H_F= \im{\, \{H^*(\phi):  H_{\mathcal E} \to H_{\mathcal F}}\}.$$
 \end{enumerate}
\end{lm}
{\em Proof.} 
In either case, $\phi^{0,q}_r=0,$ for every $q,r.$

Let us prove {\em 1.}

CLAIM.  {\em Let $r \geq 2$
 be fixed. Then
$$
{\mathcal E}_r^{p,\bullet} \lorw {\mathcal F}_r^{p, \bullet}
$$
is an   isomorphism $\forall p  \leq -r$ and  is  monic for $\forall p \leq -1.$}

The proof  is by induction on $r.$  If $r=2,$ then  the CLAIM is true by hypothesis.

Assume we have proved the CLAIM for $r.$ Let us prove it for $r+1.$

Consider the diagram
$$
\xymatrix{
    {\mathcal E}_r^{p-r,q+r-1}  \ar[r]^{d_r} \ar[d]^{\phi'_r}  &  {\mathcal E}_r^{pq} \ar[r]^{d_r} \ar[d]^{\phi_r}  & 
   { {\mathcal E}^{p+r,q-r+1}_r }  \ar[d]^{\phi''_r} \\
      {\mathcal F}_r^{p-r, q+r-1}  \ar[r]^{d_r} &   {\mathcal F}_r^{pq} \ar[r]^{d_r} & {{\mathcal F}_r^{p+r, q-r+1}}
      }
$$
Note that if $\phi'_r$ and $\phi_r$ are isomorphisms and $\phi''_r$ is monic, then
$\phi^{pq}_{r+1}$ is iso.  Since this happens by the inductive hypothesis
for $p\leq r-1,$ we have isomorphisms in that same range of $p$.

If $-r \leq p \leq -1,$ then $\phi'_r$ is an  isomorphism and $\phi_r$ is monic so that,
by Lemma \ref{pippa}, 
$\phi^{pq}_{r+1}$ is monic in this range for $p$ and the CLAIM follows.

By taking $r\gg 0,$ we get the injectivity statement for $\phi_\infty.$

Since $\phi_\infty^{0q}=0,$ $L^0 \subseteq \ke{\, H^*(\phi):H_{\mathcal E} \to H_{\mathcal F}}.$
As to the reverse inclusion we argue by contradiction. 
Let $a \in \ke{\, H^*(\phi)} \setminus L^0.$
There exists a unique $p\leq \!\!-1$ such that $0 \neq [a] \in L^pH_{\mathcal E} /L^{p+1}H_{\mathcal E}$.
One would have $\phi_\infty^{p\bullet}([a])=0,$ violating the
injectivity of the $\phi_\infty$ in this range for $p.$ This proves $1$.

The proof of  2 is analogous except, possibly, for the equality $L^1H_{\mathcal F}=\im{\,H^*(\phi})$.
Since $\phi_\infty^{0\bullet}=0,$ $L^1H_{\mathcal F} \supset \im{\,H^*(\phi})$. The reverse inclusion can be proved
as follows. Let $f \in L^1{H_{\mathcal F}}$. By the surjectivity of $\phi_\infty^{1,q},$ there exists
$e_1 \in L^1 H_{\mathcal E}$ with $(f- \phi(e_1)) \in L^2H_{\mathcal F}.$
 We repeat this procedure using the successive surjections
and conclude by induction.
\blacksquare

\medskip
Let us explain the notation in the two propositions that follow.

In the applications we have in mind, i.e. the
 proof of Theorem \ref{mtm}, the spectral sequence  ${\mathcal E}_{-i}$ of Proposition
\ref{descrfiltr}
is   a suitable realignment of the perverse spectral sequence for the cohomology group $H^*(Y_{-i}, K_{|Y_{-i}})$
of  
the element $Y_{-i}$ of a $n$-flag.    
The  intervals  $[-n,-i]$  are explained as follows.
We   work with an affine variety $Y$ of  dimension 
$n$, with a  complex $K\in {\mathcal D}_Y$,  and we choose an $n$-flag with  closed subsets
$Y_{-i}\subseteq Y$  which are affine and of dimension $n-i$. By virtue of the
Theorem \ref{cdav} on the cohomological dimension of affine varieties with respect to
perverse sheaves, the  perverse spectral sequence ${\mathcal E}_2^{s,t}$ for  the cohomology
the  complex $K_{|Y_{-i}}$
  is displayed on the
quadrants  II/III with nontrivial entries only in the columns labelled by 
 $s \in [-(n-i),0].$ The restriction maps induce maps of spectral sequences.
But this is not what we are looking for; see Remarks \ref{1122} and \ref{itsdiff}.
We want to move to the left the display of the perverse spectral sequence for $Y_{-i}$ so that
the interval of non zero $s$-columns  becomes $[-n,-i]$ 
and is thus ``aligned" with the perverse spectral sequence
for $H^*(Y, K)$, for which the corresponding interval is
$[-n,0]$. We want to do so
and still have interesting  maps of spectral sequences. This is made possible 
 by  Lemma
\ref{sotex}. The use of the re-numerated spectral sequences, i.e.
${\mathcal E}$ instead of $E$ (cf. \S\ref{ssecss}),  is not necessary,
but it makes it visually easier to analyze the realigned
maps of spectral sequences and deduce 
 strong injectivity/surjectivity properties
by virtue of the   Lefschetz hyperplane theorem  \ref{swl} and
Lemma \ref{ssis}.

In short, in Proposition \ref{descrfiltr},  the hypothesis
 (a)   mirrors  Theorem \ref{cdav} on the cohomological dimension
of affine varieties  and   the hypothesis
(b)  mirrors the  Lefschetz  hyperplane Theorem \ref{swl}.

The choice of notation in  Proposition \ref{descrfiltr2} is dictated by similar considerations,
adapted to the case of cohomology with compact supports.

\begin{pr}
\label{descrfiltr}
Let $n\geq  0$ be a fixed  integer and
$$
{\mathcal E}_{0} \lorw {\mathcal E}_{-1} \lorw {\mathcal E}_{-2} \lorw \ldots \lorw {\mathcal E}_{-n} \lorw {\mathcal E}_{-n-1}:=0
$$
be maps of spectral sequences. Denote by $\phi (i,j)= \phi_{r}^{pq}(i,j): ({\mathcal E}_{r}^{pq})_{-i} \to 
({\mathcal E}_{r}^{pq})_{-j}$ the obvious maps for $ i < j.$

Assume that 

\smallskip
\n
{\rm (a)} $({\mathcal E}^{pq}_2)_{-i}=0$ for $p \notin [-n,-i]$ and that

\smallskip
\n
{\rm (b)} 
for every $0 \leq i \leq n,$ $\phi_2^{pq} (-i, -i-1)$ is 
an iso for $p \leq -i-2$ and  monic for $p=\!-i-1$ (it is automatically zero
for $p \geq -i$ by {\rm(a)}). 

\smallskip
Then for every $\forall \; 0 \leq i \leq n:$
$$
L^{-i} H_{{\mathcal E}_{0}} = \ke{\,  \{  \phi(0, -i-1): H_{{\mathcal E}_{0}}
 \lorw H_{{\mathcal E}_{-i-1}}  \} }.
$$
\end{pr}
{\em Proof.} 
The case $n=0$ is trivial.
The system of ${\mathcal E}_{-i}$ with $i>0$ can be made to satisfy the 
hypothesis of this proposition with $n-1$ in place of $n$ by  shifting the
display one unit to the right.

The proof is by induction on $n$. 
 
The case $i=0$ is covered by Lemma \ref{ssis}, so we may assume that $i>0$.

Assume we have proved the proposition for $n-1$ and let us prove it for $n.$
The induction hypothesis translates into
$$
L^{-i} H_{{\mathcal E}_{-1}} = \ke{\,  \{ \phi(-1, -i-1): 
H_{{\mathcal E}_{-1}} \lorw H_{{\mathcal E}_{-i-1}}  \} },   \qquad \forall \; 1 \leq i \leq n.
$$
Since $L^{-i}H_{{\mathcal E}_{0}} \to L^{-i }H_{{\mathcal E}_{-1}},$ we see that
$$
L^{-i} H_{{\mathcal E}_{0}} \subseteq  \ke{ \,  \{  \phi(0, -i-1): H_{{\mathcal E}_{0}} \lorw 
H_{{\mathcal E}_{-i-1}}  \} },
 \qquad \forall \; 1 \leq i \leq n
$$
and we already know the case $i=0$.

To prove the reversed inclusion,  we argue by contradiction.
Let  $a \in \ke{\,\phi(0,-i-1)}.$ 
Note that its image $a' \in H_{{\mathcal E}_{-1}}$ maps automatically to zero into
$H_{{\mathcal E}_{-i-1}}$ so that, by the inductive hypothesis, $a' \in L^{-i}H_{{\mathcal E}_{-1}}.$

Assume  that  $a \notin L^{-i}H_{{\mathcal E}_{0}}.$
Then $\exists ! \, s>i$ such that $a \in L^{-s}H_{{\mathcal E}_{0}}  \setminus L^{-s+1} H_{{\mathcal E}_{0}},$
so that $0 \neq [a]_{-s} \in (   {\mathcal E}_{\infty}^{-s,\bullet}   )_{0} .$

By Lemma \ref{ssis}.1, the image  of this element, $[a']_q $ in $({\mathcal E}_{\infty}^{-s, \bullet})_{-1}$ is not zero.
This would imply that $ a' \notin L^{-s+1} H_{{\mathcal E}_{-1}} \supseteq
L^{-i} H_{{\mathcal E}_{-1}},$ a contradiction.
\blacksquare

\medskip
Analogously, we have the following 

\begin{pr}
\label{descrfiltr2}
Let $n\geq  0$ be a fixed  integer and
$$
0 =: {\mathcal E}_{n+1} \lorw {\mathcal E}_{n} \lorw {\mathcal E}_{n-1} \lorw \ldots \lorw {\mathcal E}_{1} \lorw {\mathcal E}_{0}
$$
be maps of spectral sequences. Denote by $\phi (i,j): {\mathcal E}_{i} \to {\mathcal E}_{j}$ the obvious maps
for $ i\geq j.$

Assume that 

\smallskip
\n
{\rm (a)}
$({\mathcal E}^{pq}_2)_{i}=0$ for $p \notin [i,n]$ and that

\smallskip
\n
{\rm (b)} for every $0 \leq i \leq n,$ $\phi_2^{pq} (i, i-1)$ is 
an iso for $p \geq i+1$ and  epic  for $p=i$ (it is automatically zero
for $p \leq i-1$ by {\rm (a)}).

\smallskip
Then 
$$
L^{i} H_{{\mathcal E}_{0}} = \im{ \, \{ \phi(i, 0): H_{{\mathcal E}_{i}} \lorw H_{{\mathcal E}_{0}}  \} }, 
\qquad \forall \; 0 \leq i \leq n.
$$
\end{pr}
{\em Proof.} Completely analogous to the proof of Proposition
\ref{descrfiltr}, via  Lemma \ref{ssis}.2.
 \blacksquare

\subsubsection{Jouanolou Trick: reduction to the affine case }
\label{ssjo}
In this section, we prove Lemma \ref{pjc} which   is one 
way to  reduce the study of the perverse spectral sequence
on a quasi projective variety $Y$, to the case of affine varieties
by replacing  $Y$ with the affine $\mathcal Y$.
For a different approach, using two flags,  see \ci{decmigpf}.

For the  notions of translated spectral sequence ${\mathcal E}(l)$ 
with   filtration  
$L(l)$ see  \S\ref{ssecss}.

Let $Y$ be a quasi projective variety of dimension $n$, let
$f : X \to Y$ be a map and let
$K \in {\mathcal D}_Y,$ $C \in D_X$. 

By  fixing  an ${\Bbb A}^d$-fibration $\pi: {\mathcal Y} \to Y$
as in Proposition \ref{jo}, we obtain  the
Cartesian diagram
\begin{equation}
\label{cdjo}
\xymatrix{
{\mathcal X} \ar[r]^\pi \ar[d]^f &  X \ar[d]^f \\
{\mathcal Y} \ar[r]^\pi & Y.
}
\end{equation}
Since the fibers of the maps $\pi$ are affine spaces, we have canonical identifications
\begin{equation}
\label{pppp}
Id \simeq \pi_{*} \pi{^{*}}, \quad 
H^*(Y, K) = H^*({\mathcal Y}, \pi^* K),
\qquad 
\pi_{!} \pi^{!} \simeq Id, \quad H^*_c({\mathcal Y}, \pi^! K) = H^*_c(Y, K);
\end{equation}
in fact,  the first one follows from \ci{ks}, Corollary 2.7.7.(ii), the second follows from the first one, the fourth from the third and the third from
the first one in view of the fact that  Poincar\'e-Verdier duality  exchanges the pull-back
$\pi^*$  with the extraordinary pull-back $\pi^!$.

Since $\pi$ is smooth of relative dimension $d,$ we have
a canonical identification of  functors
\begin{equation}
\label{tuko}
\pi^* [d] = \pi^! [-d].
\end{equation}
The functors $\pi^* [d] = \pi^![-d]$  
are $t$-exact with respect to middle perversity. The functors
$\pi^*= \pi^! [-2d]$ are exact in the  usual sense.

Given a flag ${\mathcal Y}_*$ on ${\mathcal Y}$, we have the pre-image flag
${\mathcal X}_* := f^{-1} {\mathcal Y}_*$ on ${\mathcal X}$. Recall that we have the two associated filtrations $F_{{\mathcal Y}_*}$ and $G_{ {\mathcal Y}_*}$ of $\S$\ref{ssflag}. 

\begin{lm}
\label{pjc}
There are canonical identifications  of  filtered  abelian groups:
\begin{equation}
\label{ffqq1}
H^*(Y, K)  =
H^*({\mathcal Y}, \pi^{*}K), \qquad    L^Y_{\pp} =L^{\mathcal Y}_{\pp }(-d) ,
\qquad  L^Y_\tau = L^{\mathcal Y}_\tau,
\end{equation}
\begin{equation}
\label{ffqq2}
H_{c}({\mathcal Y}, \pi^{!}K) =
H_{c}( Y , K), \qquad          L^Y_{\pp } = L^{\mathcal Y}_{  \pp } (d),
\qquad  L^Y_\tau = L^{\mathcal Y}_\tau (2d).
\end{equation}
We have the following relations in $H^*(X, C)$:
\begin{equation}
\label{ffqq3}
F_{{\mathcal X}_*} \supseteq F_{{\mathcal Y}_*}, \qquad 
L^{f: X \to Y}_{\pp} = L^{f: {\mathcal X} \to {\mathcal Y} }_{\pp} (-d), \qquad
L^{f: X \to Y}_\tau = L^{f : {\mathcal X} \to {\mathcal Y} }_\tau,
\end{equation}
and  the following relations in $H^*_c(X, C):$
\begin{equation}
\label{ffqq4}
G_{{\mathcal X}_*} \subseteq  G_{{\mathcal Y}_*},
\qquad 
L^{f: {\mathcal X} \to {\mathcal Y} }_{\pp} (d) =
L^{f: X \to Y}_{\pp} ,
\qquad
L^{f: X \to Y}_{\tau}= L^{f: {\mathcal X} \to {\mathcal Y} }_{\tau} (2d).
\end{equation}
\end{lm}
{\em Proof.}
The statements (\ref{ffqq1}) and (\ref{ffqq2})
follow from (\ref{pppp}), the $t$-exactness
of  $\pi^*[d] = \pi^! [-d]$,  the exactness of $\pi^* = \pi^![-2d]$ and from
Lemma \ref{sotex}.

The inclusion  in (\ref{ffqq3}) is seen as follows.
There is the commutative diagram
\begin{equation}
\label{ttt}
\xymatrix{
H^*( {\mathcal X} , \pi^* C) \ar[d]^{a}      & = &  &  H^*(   {\mathcal Y}, f_* 
\pi^* C) \ar[d]^{a'} \\
H^*( {\mathcal X}_{s}, i_s^* \pi^* C  )  & = &  H^*({\mathcal Y}_s, f_* i_s^* \pi^* C) & 
H ( {\mathcal Y}_s, i_s^* f_* \pi^* C) \ar[l]_{b}  ,
}
\end{equation}
where $b$ stems from the base change map (\ref{bc1})
$i_s^* f_* \to f_* i_s^*.$
The kernels
of the vertical restriction  maps $a$
and $a'$ define the filtrations $F_{\mathcal X}$ and $F_{\mathcal Y}$
and it is clear that $\ke{\, a} \supseteq \ke{\, a'}.$

The second (third, respectively) equality in (\ref{ffqq3}) follows from the definition
of $L^f_{\pp}$ (of $L^f_\tau,$ respectively)  the smooth base change isomorphism
 $f_* \pi^* = \pi^* f_*$
for the smooth map $\pi$ and the second (third, respectively) equality
 in (\ref{ffqq1}).

The proof of (\ref{ffqq4}) runs parallel to the one just given
via the use of the base change map  $f_! i_s^! \to i_s^! f_!$
and the base 
change isomorphism
$f_! \pi^! = \pi^! f_!.$
The  inclusion is  reversed with respect  
to the one in cohomology, and this is because   
the flag filtration in cohomology with compact supports is, by definition, given by the images
of the co-restriction maps. 
\blacksquare

\begin{rmk}
\label{atpmc}
{\rm
The proof of Lemma \ref{pjc}
makes it clear that the failure of the base change map to be an isomorphism
is responsible for the inequality
$F_{{\mathcal X}_*} \neq F_{{\mathcal Y}_*}.$
If $f: X \to Y$ is proper,  then, the base change isomorphism
for proper maps yields the equality  $F_{ {\mathcal X}_*} = F_{ {\mathcal Y}_* }$ 
in cohomology, as well as in cohomology with compact supports.
  If $f$ is not proper, then one can still have
the equality $F_{ {\mathcal X}_*} = F_{ {\mathcal Y}_* }$ in cohomology, provided 
that, for every element $i_s: {\mathcal Y}_s \to {\mathcal Y}$ of the flag, the
base change map
$i_s^* f_* \to f_* i_s^*$ is an isomorphism when evaluated on
$\pi^*C.$  Similarly, in the case of cohomology with compact supports. As we show in 
\ci{decbday},  and also in \ci{decmigpf}, this can be achieved by choosing a system
of general hyperplane sections. In particular, the following holds.
}
\end{rmk}
\begin{pr}
\label{1ar}
In the situation of {\rm Lemma \ref{pjc}}, if the flag ${\mathcal Y}_*$ is chosen to be general, then we have the equality $F_{ {\mathcal X}_*} = F_{ {\mathcal Y}_* }$.
\end{pr}

\subsection{The geometry of the perverse and perverse Leray  filtrations}
\label{tgofsss}

Let $Y$ be a quasi projective variety, let  $K \in {\mathcal D}_Y$,
let $f: X \to Y$ be a map, and let  $C \in {\mathcal D}_X$.

Recall that, the perverse
Leray filtration on $H^*(X,C)$ is defined to be the perverse
 filtration on $H^*(Y,f_*C)$. Similarly, for $H_c^*(X,C)=H_c^*(Y, f_!C)$.
 
In this  section, we employ  the set-up of $\S$\ref{ssjo}
and 
we identify:  
\begin{itemize}
\item the perverse filtrations on $H^*(Y, K)$ and on  $H^*_c(Y,K)$
with suitable flag filtrations on the  auxiliary  affine  variety ${\mathcal Y}$,
 and  
 \item
  the  perverse Leray filtrations  on $H^*(X, C)$ and  on $H^*_c(X,C)$
with suitable flag filtrations on the  auxiliary   variety  ${\mathcal X}$.
\end{itemize}

Lemma \ref{pjc}, relates the perverse Leray filtration on $H^*(X, C)$  with the one
on $H^*({\mathcal X},\pi^*C)$ and similarly for compact supports. 
We employ the following identifications
\[ H^*(X,C) = H^*(Y, f_* C) = H^*( {\mathcal Y}, \pi^* f_* C = f_* \pi^* C) = H^*({\mathcal X}, \pi^* C),
\]
\[ H^*_c(X,C) = H^*_c(Y, f_! C) = H^*_c( {\mathcal Y}, \pi^! f_! C = f_! \pi^! C) = H^*({\mathcal X}, \pi^! C),
\]
of cohomology groups and the following ones for the  filtrations on them
\begin{equation}
\label{zz1}
L^{f:X\to Y}_{\pp} := L^Y_{\pp} \stackrel{(\ref{ffqq1})}= L^{\mathcal Y}_{\pp} (-d) = : L^{f: {\mathcal X} \to {\mathcal Y}}_{\pp} (-d), \quad
\mbox{on} \; H^*(X, C),
\end{equation}
\begin{equation}
\label{zz2}
L^{f:X\to Y}_{\pp} := L^Y_{\pp} \stackrel{(\ref{ffqq2})}= L^{\mathcal Y}_{\pp} (d) = : L^{f: {\mathcal X} \to {\mathcal Y}}_{\pp} (d), \quad
\mbox{on} \; H^*_c(X, C).
\end{equation}

\subsubsection{The perverse  filtrations on $H^*(Y,K)$ and $H^*_c(Y,K)$}
\label{tpfojimo}

\begin{tm}
\label{mtm} Let $Y$ be a quasi projective variety of dimension $n$ and $K \in {\mathcal D}_Y$.
Let $\pi: {\mathcal Y} \to Y$ be  an ${\Bbb A}^d$-fibration with  ${\mathcal Y}$ affine. 
There exists a $(n+d)$-flag ${\mathcal Y}_*$
on ${\mathcal Y}$ such that
we have  an equality of filtered abelian groups
\begin{equation}
\label{mtmco}
 H^*(Y, K) = H^*({\mathcal Y}, \pi^* K), \qquad  
L^Y_{\pp} = F_{ {\mathcal Y}_* } (-d),
\end{equation}
\begin{equation}
\label{cptsupmtm}
 H_c^* (  Y,  K)= H^*_c({\mathcal Y}, \pi^! K) , \qquad
L^Y_{\pp} = G_{ {\mathcal Y}_* } (d).
\end{equation}
If $Y$ is  an affine variety, then  one may take ${\mathcal Y} = Y$ and $d=0$, and then
\begin{equation}
\label{mtmcozz}
   L^Y_{\pp} = F_{ Y_* }  \;\; \mbox{on} \;\;
H^*(Y, K), \qquad 
L^Y_{\pp} = G_{ Y_* }  \;\; \mbox{on} \;\;
H^*_c(Y, K).
\end{equation}
\end{tm}
{\em Proof.}
The proofs for cohomology and for cohomology with compact supports run parallel.
By virtue of  the equalities
${L^{Y}_{\pp}} = L^{ {\mathcal Y}}_{\pp}(-d)$  (\ref{ffqq1}) 
and
${L^{Y}_{\pp}} = L^{ {\mathcal Y}}_{\pp}(d)$  (\ref{ffqq2}) 
in Lemma \ref{pjc}, we
are left with showing that we can choose ${\mathcal Y}_*$ on ${\mathcal Y}$ so that
${L^{\mathcal Y}_{\pp}} = F_{{\mathcal Y}_*}$ in cohomology,
and 
${L^{\mathcal Y}_{\pp}} = G_{{\mathcal Y}_*}$ in cohomology with compact supports.
In particular, we may assume that ${\mathcal Y} =Y$ is affine.

Let  $\Sigma$ be a stratification of $Y$ with the property that
$K$ is $\Sigma$-constructible, i.e.
$K \in {\mathcal D}^{\Sigma}_{Y}.$

Let $i: Y_{-1}  \to Y$ be a general hyperplane section of $Y.$
Then, for every $m \in \zed,$  
the complexes $i^{*} \pc{m}{K}[-1] = i^{!} \pc{m}{K}[1] 
\in {\mathcal P}_{ Y_{-1}},$ i.e. they
are perverse on $Y_{-1}.$

By  Theorem \ref{cdav} on the cohomological dimension
of affine varieties   and by the  Lefschetz  hyperplane Theorem \ref{swl},
 there is a general (cf. Remark \ref{usofrutto}) hyperplane section $i:=i_{-1}: Y_{-1} \to Y$
 such that, for every $m \in \zed$, the natural maps 
 $$
 \xymatrix{
 a^{j}: H^{j}(Y, \pc{m}{K}) \ar[r] & H^{j}(Y_{-1}, i^{*} \pc{m}{K}), 
 \\
 b^{j}: H^{l}_{c}(Y_{-1}, i^{!} \pc{m}{K}) \ar[r] & H^{l}_{c}(Y, \pc{m}{K}),
 }
 $$
satisfy the following conditions:
\begin{equation}
\label{es1}
H^{j}(Y,\pc{m}{K})= 0, \quad j \notin [-n,0], \qquad 
H^{j}(Y_{-1}, i^{*}\pc{m}{K})= 0, \quad j \notin [-n,-1],
\end{equation}
\begin{equation}
\label{es2}
\mbox{    $a^{j}$ is iso for $j \in [-n,-2]$}, \qquad 
\mbox{   $a^{-1}$ is monic};
\end{equation}
\begin{equation}
\label{es1c}
H^{j}_{c}(Y,\pc{m}{K})= 0, \quad j \notin [0,n], \qquad 
H^{j}_{c}(Y, i^{!}\pc{m}{K})= 0, \quad j \notin [1,n],
\end{equation}
\begin{equation}
\label{es2c}
\mbox{    $b^{j}$ is iso for $j \in [2,n]$}, \qquad 
\mbox{   $b^{1}$ is epic}.
\end{equation}
Note that the complex $i^{*}K$ is constructible with respect 
to the stratification $\Sigma_{-1}$ on $Y_{-1}$ induced by $\Sigma$
and that $i^{*}[-1]= i^{!}[1] : {\mathcal D}^{\Sigma}_{Y} \to {\mathcal D}^{ \Sigma_{-1}}_{Y_{-1}}$ 
are $t$-exact.

We iterate the construction and take $Y_{-l-1}$
to be a general hyperplane section of $Y_{-l}.$

Let ${\mathcal E}_{-l}$ be the perverse spectral sequence (\ref{per})  for
$H^*(Y_{-l}, K_{| Y_{-l}})$ translated   by $+l: $ 
\begin{equation}
\label{tadpp}
{\mathcal E}_{-l} := {\mathcal E} ( K_{| Y_{-l}},    \pp  ) (+l). 
\end{equation}
 
Let $_{c}{\mathcal E}_{l}$ be the perverse spectral sequence (\ref{per})  for
$H_{c}(Y_{-l}, i^{!}_{-l}K)$ translated   by $-l: $ 
$$
_{c}{\mathcal E}_{l} :=  \, _{c}{\mathcal E} ( i^{!}_{-l}K)  ,    \pp  ) (-l). 
$$

Since  $i^*_{-l}[-l]= i^{!}_{-l}[l]$ are $t$-exact, 
 Lemma
\ref{sotex} yields two  systems of maps of spectral sequences
$$
{\mathcal E}_{0} \lorw {\mathcal E}_{-1} \lorw  \ldots \lorw {\mathcal E}_{-n} \lorw {\mathcal E}_{-n-1}=0,
$$
$$
0= {\mathcal E}_{n+1} \lorw {\mathcal E}_{n} \lorw  \ldots  \lorw {\mathcal E}_{1} \lorw {\mathcal E}_{0}.
$$
By (\ref{es1}) and (\ref{es2}),  
this system of spectral sequences satisfies the hypotheses of 
Proposition \ref{descrfiltr} and
the conclusion for cohomology  follows.

By (\ref{es1c}) and (\ref{es2c}),  
this system of spectral sequences satisfies the hypotheses of 
Proposition \ref{descrfiltr2} and
the conclusion for cohomology with compact supports   follows. 
\blacksquare

\begin{rmk}
\label{efmtm}
{\rm 
Let us write out the conclusions above when $Y$ is affine:
$$
L^{p}_{\pp} \;H^{j}(Y,K)=
\ke{\,
\{
H^{j}(Y,K) \lorw
H^{j}(Y_{p -1}, K_{| Y_{p-1}  }
\}
}, \qquad \forall \, p \in \zed,
$$
$$
L^{p}_{\pp} \; H^{j}_{c}(Y, K) = \im{\,
\{
H^{j}_{c}(Y_{p}, i_{p}^{!} K) \lorw
H^{j}_{c}(Y, K)
\},
} \qquad \forall \, p \in \zed.
$$
}
\end{rmk}

\begin{rmk}
\label{pul}
{\rm
As the proof shows, one can choose the flag ${\mathcal Y}_*$
by taking $(n+d)$ general hyperplane sections of ${\mathcal Y}$. 
In fact, one can   take $(n+d)$ general hypersurface
sections of ${\mathcal Y}$ of varying degrees.
The subspaces of the  filtrations can be viewed, in cohomology,  as kernels of the restriction maps and,
in cohomology with compact supports, as  images 
of the  co-restriction maps  associated with
the maps
 $\pi_l:= \pi \circ i_l: {\mathcal Y}_l
 \to Y$. 
}
\end{rmk}

\begin{rmk}
\label{loffi}
{\rm ({\bf Type of perverse filtrations})
By trivial reasons of indexing,  the flag filtration $F_{Y_*}$
is of type $[-n,0]$ on each group $H^*(Y,K)$, and
Theorem \ref{cdav} on the cohomological dimension of affine varieties implies that
if $Y$ is affine, then the perverse filtration $L_{\pp}$ on the group
$H^*(Y,K)$  is
of type $[-n,0]$.  This is in accordance with Theorem \ref{mtm}, i.e.
if  $Y$ is affine of dimension $n$  and $Y_*\subseteq Y$ is  a  general
$n$-flag of linear sections, then
$L_{\pp} =F_{Y_*}$. Similarly for   cohomology with compact supports,
where the type   of $G_{Y_*}$ on the affine $Y$ is $[0,n]$.
If $Y$ is quasi projective,  then Theorem \ref{mtm} yields
    $[-n, d]$ and $[-d,n]$ as bounds on the type  of the perverse filtrations
    $L_{\pp}$ in cohomology and in cohomology with compact supports,
    respectively, where $d$ is the dimension
 of the fiber ${\Bbb A}^d$ of the chosen  map $\pi : {\mathcal Y} \to Y$.
These bounds on $L^Y_{\pp}$ are 
{\em not } sharp. One needs to replace ${d}$ by $t,$
where $t+1$ is the smallest  number of affine open sets 
for an   open  covering of $Y$
given by affine open subvarieties. 
However, this does not  create problems in  applications.
}
\end{rmk}

\subsubsection{The perverse Leray filtrations on  $H^*(X,C)$ and $H^*_c(X,C)$}
\label{sstplfhhc}

We are aiming at a geometric description of the perverse Leray filtrations using
a flag  on ${\mathcal X}$. The obvious candidate is the flag ${\mathcal X}_*:= f^{-1} {\mathcal Y}_*$.

Recall that the  inclusions (\ref{ffqq3}) and (\ref{ffqq4}) in Lemma \ref{pjc}
can be  strict in view of the possible failure of the relevant base change
of a general flag ${\mathcal Y}_*$ on ${\mathcal Y}$ corrects this failure.

\begin{tm}
\label{flplr}
Let $f: X \to Y$ be a map of varieties where $Y$ is quasi projective of dimension $n$.
Let $\pi: {\mathcal Y} \to Y$ be an ${\Bbb A}^d$-fibration with  ${\mathcal Y}$ affine
and 
let ${\mathcal X}= {\mathcal Y} \times_Y X$.
Then there is a $(n+d)$-flag ${{\mathcal X}_{*}}$ on $\mathcal X$
such that
$$
  L^{f:X \to Y}_{\pp}  =  F_{  {\mathcal X}_{*}  } (-d)  
  \;\; \mbox{on $H^*(X, C)$},
   \quad  \mbox{and} \quad 
  L^{f: X \to Y}_{\pp} = G_{  {\mathcal X}_{*}  }     (d)  \;\;
\mbox{on $H^*_c(X, C).$} 
 $$
If $Y$ is an  affine variety, then we may take ${\mathcal X}= X$  and $d=0$ and then
$$   
  L^{f: X \to Y}_{\pp}  =  F_{  { X}_{*}  }  
  \;\; \mbox{on $H^*(X, C)$},
   \quad  \mbox{and} \quad 
L^{f:X \to Y}_{\pp}  = G_{  X_{*}  }       \;\;
\mbox{on $H^*_c(X, C).$} 
 $$
\end{tm}
{\em Proof.} We prove the result for cohomology.
The case of compact supports is analogous. 
Choose a $(n+d)$-flag ${\mathcal Y}_*$ on ${\mathcal Y}$ using
$(n+d)$ general hyperplane sections of ${\mathcal Y}$ such that,
 having set ${\mathcal X}_* := f^{-1} {\mathcal Y}_*$,
we have:
\[ L^{\mathcal Y}_{\pp}\stackrel{Thm. \ref{mtm}}= F_{ {\mathcal Y}_* }
\stackrel{Prop. \ref{1ar}}= F_{ {\mathcal X}_*}.\]
We conclude by  using  the equality  (\ref{zz1}) $L^{\mathcal Y}_{\pp}(-d) = L^{f:X \to Y}_{\pp}$.
\blacksquare

\begin{rmk}
\label{tyyt}
{\rm 
If  $Y$ is affine, then we have
 $$
 ( L^{f}_{\pp}  )^{p } H^{j}(X, C) = \ke{ \,  \{   
 H^*(X, C) \lorw H^{j}(X_{p-1}, C_{|  X_{p-1}}) \},
 }
 $$
  $$
 ( L^{f}_{\pp}  )^{p } H^{j}_{c}(X, C) = \im{ \, \{   
 H^{j}_{c}(X_{p}, i^{!}_{p}C) \lorw H^{j}_{c}(X, C) \}.
 }
 $$
}
\end{rmk}

\begin{rmk}
\label{noqp}
{\rm
Note that in Theorem \ref{flplr} we do not need to assume that
$X$ is quasi projective, nor that $f$ is proper.
}
\end{rmk}

\begin{rmk}
\label{ar}
{\rm
D. Arapura \ci{arapura} proved that if $f: X \to Y$ is a projective map of quasi projective varieties,
then 
the  {\em standard} Leray filtration on the cohomology  groups $H^*(X,\zed)$  can be described geometrically
using flags in {\em special} position. His methods  are different
from the ones  of the present paper (and of \ci{decmigpf}). However, I feel a
great intellectual debt to \ci{arapura}.
I do not  not know of a way to describe the standard Leray filtration in cohomology
using flags on (varieties associated with) the domain of a non proper map.
The paper \ci{decbday} is devoted to provide such a description
for the Leray filtration on cohomology with compact supports
via compactifications of varieties and of maps.
}
\end{rmk}

\section{Applications of the results on filtrations}
\label{secappl}
Recall that  the acronym MH(S)S stands for mixed Hodge (sub)structure.
Due to the functorial nature of the canonical MHS on algebraic varieties,
the description of  perverse filtrations in terms of flags, i.e.
as kernel of restrictions and as images of co-restrictions
to subvarieties, is amenable to applications to the mixed Hodge theory of the cohomology
and intersection cohomology of 
quasi projective varieties. In this section we work out some of these applications.

\subsection{Perverse Leray and MHS: singular cohomology}
\label{mhsapl}
The following  theorem appears in a stronger form
(involving filtered complexes and spectral sequences)
 in \ci{decmigpf} and we include it here  for the reader's  convenience
as  it is an important preliminary result to the applications that follow.

\begin{tm}
\label{hmhspl}
Let $f:X \to Y$ be a map of algebraic varieties
 with $Y$ quasi-projective.

\n
There are an integer $d$, a   variety ${\mathcal X}$, and   a flag ${\mathcal X}_*$
on ${\mathcal X}$
 such that there are identities of   filtered groups
$$
\xymatrix{
(H^*(X, \zed), L_{\pp}^f  )  &=&  ( H^*({\mathcal X}, \zed), F_{ {\mathcal X}_* }(-d) ),\\
(H^*_c(X, \zed), L_{\pp}^f  ) &=&  ( H^*_c({\mathcal X}, \zed), G_{ {\mathcal X}_* }(d) ).
}
$$
In particular,  the perverse Leray filtrations on $H^*(X,\zed)$ and $H^*_c(X, \zed)$ are  by MHSS.
\end{tm}
{\em Proof.} The first statement is a mere application of 
Theorem \ref{flplr} to the 
case $C= \zed_X$, and  
the subspaces
of the  perverse filtrations are   the kernels of  the restriction maps 
$H^*(X,\zed) = H^*({\mathcal X}, \zed) \to H^*({\mathcal X}_p, \zed)$ and the  images of 
the co-restriction maps $H^{*+ 2d+2p}_c({\mathcal X}_p,  \zed) \to H^{*+2d}_c({\mathcal X},  \zed ) =
H_c^{*}(X,\zed)$,
respectively.  
The second statement about MHSS  follows
from the usual functoriality properties of the MHS on the cohomology groups
and  on the cohomology groups  with compact supports of  varieties \ci{ho3}.
\blacksquare

\begin{rmk}
\label{lm11}
{\rm 
The case of the Leray filtration in cohomology    is dealt with in 
\ci{decbday} where the relation with Arapura's results
\ci{arapura} is discussed.
}
\end{rmk}

\begin{rmk}
\label{dtosiato}
{\rm 
Theorem \ref{hmhspl} holds if we replace cohomology with intersection cohomology; see Proposition \ref{cor-555}. The proof is  formally
analogous.
However,  we must first 
endow   intersection cohomology with a MHS (Theorem 
 \ref{tm111})
and then verify   the compatibility of the  resulting MHS
with restrictions to general hyperplane sections
(see the proof of Proposition \ref{cor-555}).  
}
\end{rmk}

\subsection{Review of the decomposition theorem} 
\label{dtmhs}
In the sequel of this paper, we employ $\rat$-coefficients (we still denote the corresponding
category ${\mathcal D}_Y$). 
The reason for employing rational coefficients stems from the use of the decomposition
theorem, which does not hold over the integers. Moreover, we use Poincar\'e-Verdier
duality  in the form of perfect pairings between rational vector spaces.

We also  assume that $X$ and $Y$
are irreducible and quasi projective of dimension
$m$ and $n$, respectively. This  only  because it makes the exposition simpler.
All the results we prove hold without the irreducibility assumption, 
with pretty much the same proofs.
 
Let $f: X \to Y$ be a proper map. 
The decomposition theorem, due to Beilinson-Bernstein-Deligne-Gabber  (see the survey \ci{decmigbams})
yields the existence of   direct sum decompositions in ${\mathcal D}_Y$ for the direct image of the
(rational) intersection complex of  $X$:
\begin{equation}
\label{dtgl}
\phi: \;  f_* IC_X \; \simeq  \;
\bigoplus_{i,l,S}{ IC_{\overline{S}} (L_{i,l,S}) [-i]}, 
\end{equation}
where $i \in \zed$, $ 0 \leq  l \leq n$,   $Y = \coprod_{ 0 \leq l \leq n} S_l$ is a stratification
of $Y$, part of a stratification for the map $f$,
 $S_l$ is  the non-necessarily connected $l$-dimensional stratum,
$S$ ranges over the set of connected components of  $S_l$, and
for every $i, l, S$,  the symbol
$L_{i,l,S}$ denotes a  semisimple local system of rational vector spaces on $S$.

By grouping the summands  with the same shift $[-i]$, we obtain  canonical identifications
\begin{equation}
\label{grouper}
\pc{i}{ f_* IC_X } = \bigoplus_{l,S}{ IC_{\overline{S}} (L_{i,l,S})}.
\end{equation}

Recall that $I\!H^j(X) = H^{j-m}(X, IC_X)$ and that
$I\!H^j_c(X) = H^{j-m}_c(X, IC_X)$.

We define the perverse cohomology groups of $X$ (with respect to $f$) by setting
\begin{equation}
\label{0011}
I\!H^j_i(X) :=  H^{j-i-m} ( Y, \pc{i}{ f_* IC_X } ),
\qquad
I\!H^j_{i,l,S} (X): = \bigoplus_{l,S} H^{j-i-m} (Y, IC_{\overline{S}} (L_{i,l,S}) ). 
\end{equation}
While there may be no natural choice for the isomorphism
$\phi $ in (\ref{dtgl}), the perverse cohomology groups are natural
subquotients of the  groups $I\!H^*(X)$.

We define the perverse cohomology groups with compact supports in a similar way.

We have   canonical identifications
\begin{equation}
\label{mnb}
I\!H^j_i (X) = \bigoplus_{l,S} I\!H^j_{i,l,S}(X), \qquad
I\! H^j_{c, i} (X) = \bigoplus_{l,S} I\!H^j_{c, i,l,S}(X).
\end{equation}

Recall that if $X$ is  nonsingular, then $IC_X = \rat_X [m]$
and $I\!H^j(X) = H^j(X)$, etc.  In this case,
we denote the perverse cohomology groups as follows: 
\begin{equation}
\label{-0-0}
H^j_i(X), \qquad H^j_{c,i} (X), \qquad H^j_{i,l,S}(X), \qquad
H^j_{c, i,l,S}(X).
\end{equation}

\subsection{Decomposition theorem and mixed Hodge structures}
\label{dtamhs}

In this section we state   some  mixed-Hodge-theoretic results
concerning the  intersection cohomology of quasi projective varieties.
These results are proved in $\S$\ref{schema} and $\S$\ref{-1} by making use of the geometric description of the perverse Leray filtration given in Theorem \ref{hmhspl}.

\begin{tm}
\label{tm111}
{\rm {\bf (MHS on intersection cohomology)}}
Let $Y$ be an irreducible  quasi projective variety of dimension $n$. Then
\begin{enumerate}
\item
The groups $I\!H^*(Y)$ and $I\!H^*_c(Y)$ carry a canonical  MHS.
If $Y$ is nonsingular, then this MHS coincides with
Deligne's \ci{ho3}. If $f: X \to Y$ is a resolution of the singularities
of $Y$, then the MHS on $I\!H^*(Y)$ and 
$I\!H_c^*(Y)$  are  canonical subquotient MHS of the MHS on $H^*(X)$ and $H^*_c(X)$,
respectively.
\item
The Goresky-MacPherson-Poincar\'e duality isomorphism yields  isomorphism of MHS
\[I\!H^{j}(Y) \simeq I\!H^{2n-j}_c(Y)^{\vee} (-n).\]  
\item
The canonical maps $a: H^*(Y) \to I\!H^*(Y)$,
$a': H_c^*(Y) \to I\!H_c^*(Y)$
 are  maps of MHS.
If $Y$ is projective, then  $\ke \, \{a: H^j(Y) \to I\!H^j(Y) \} =
W_{j-1} H^j(Y)$ (the  subspace of weights $\leq j-1$).
\end{enumerate}
\end{tm}

\begin{tm}
\label{ffgg}
{\rm {\bf (Pieces of the decomposition theorem and MHS)}}
Let $f: X \to Y$ be a projective map of quasi projective irreducible  varieties
and let $m:= \dim{X}$.
Then:
\begin{enumerate}
\item
The subspaces of the perverse Leray filtrations
on $I\!H^*(X)$ and on $I\!H^*_c(X)$ are MHSS for the MHS of 
Theorem {\rm \ref{tm111}}.
\item
The perverse cohomology groups  $I\!H^j_i(X)$ and $I\!H_{c,i}^j(X)$  in (\ref{0011})
carry natural MHS which are subquotients of the natural  MHS 
on $I\!H^j(X)$ and on $I\!H^j_c(X)$, respectively.
\item
The subspaces 
 $I\!H^j_{i,l,S}(X) \subseteq I\!H^j_i(X)$ and 
 $I\!H^j_{c,i,lS}(X) \subseteq I\!H^j_{c,i}(X)$  in  (\ref{mnb}) are MHSS.
 \item 
The Poincar\'e pairing  isomorphisms of Theorem {\rm \ref{tm111}} applied to $X$ descend
to the perverse cohomology groups and induce 
    isomorphisms of MHS
\[ I\!H^j_{i}(X) \simeq  I\!H^{2m-j}_{c,-i}  (X)^{\vee} (-m),
\qquad
I\!H^j_{i,l,S}(X) \simeq  I\!H^{2m-j}_{c,-i,l,S}  (X)^{\vee} (-m).
\]
 \end{enumerate}
 \end{tm}

\begin{tm}
\label{appl3}
{\rm {\bf (Hodge Theoretic splitting $\phi$) }}
There exist splittings $\phi$ in {\rm (\ref{dtgl})} for which the splittings
$$
\phi\, :\,  I\!H^j(X, \rat) \simeq  \bigoplus_{i,l,S}{ \, I\!H^j_{i,l,S}(X) },
\qquad \qquad \phi \, : \, I\!H^j_c(X, \rat) \simeq \bigoplus_{c,i,l,S}{ 
 I\!H^j_{c,i,l,S}(X) }
$$
are isomorphisms of MHS for the MHS in Theorem {\rm \ref{tm111}}, part {\rm 1},  and 
Theorem {\rm \ref{ffgg}}, part {\rm 3}.
\end{tm}

\begin{tm}
\label{appl5}
{\rm {\bf (Induced morphisms in intersection cohomology)}}
Let $f: X \to Y$ be a proper map
of quasi projective irreducible varieties.  There is a 
canonical splitting injection
$$
\gamma: IC_Y \lorw \pc{\dim{Y} - \dim{X} }{f_* IC_X}
$$
and there is a  choice of $\phi$ in {\rm (\ref{dtgl})} that yields a commutative diagram
of MHS:
$$
\xymatrix{
H^j(Y) \ar[r]^{a_Y} \ar[d]^{f^*} & I\!H^j (Y) \ar[d]^{ \phi \gamma   } \\
H^j(X) \ar[r]^{a_X} & I\! H^j (X).
}
$$
\end{tm}

Saito's work \ci{samhm,saitomathann}  on mixed Hodge modules
(MHM)  implies all the 
results on MHS we prove in $\S$\ref{secappl}, with one caveat:
it is not a priori clear that the MHS coming from the theory of  
MHM coincide with the ones of this paper.

\begin{tm}
\label{coincide}
{\rm {\bf (Comparison with M. Saito's MHS)}}
The MHS appearing in Theorems \ref{tm111} and \ref{ffgg}
coincide with the ones arising from MHM.
\end{tm}

\subsection{Scheme of proof of the results in $\S$\ref{dtamhs}}
\label{schema}
With the exception of Theorem \ref{coincide}, all the results
listed in $\S$\ref{dtamhs} were proved in \ci{decmightam, decmigseattle}
in the case when $X$ and $Y$  are  projective varieties
(in which case all the  Hodge structures in question are pure).

Roughly speaking, we first prove the results in \ref{dtamhs}
in the special case when $X$ is nonsingular, and then we use resolution
of singularities to conclude. 

Let me outline more precisely the structure of the proofs.  In fact, in what follows 
we prove these results, except for certain assertions which are then  proved
in $\S$\ref{-1}, Propositions \ref{pfbmhs}, \ref{cor-555} and  Lemmata \ref{indst}, \ref{crucial}.  Many of the details carry over verbatim from
the projective case and will not be repeated here; we simply point the reader to the original proofs
in \ci{decmightam,decmigseattle}.  Some other details which seem to be  less routine are spelled-out.

\begin{enumerate}
\item
We prove 
Theorem \ref{ffgg}, parts 1, 2 and 4  in the  case when $X$ is nonsingular by
using the geometric description of the perverse Leray filtration given by
Theorem \ref{hmhspl}. This is done in  Proposition \ref{pfbmhs}.

\item
We prove Theorem \ref{ffgg} part 3 in the case when $X$ is singular
in  Lemmata \ref{indst} and \ref{crucial}. These  two lemmata  are adapted
from \ci{decmightam}, Lemma 7.1.1  and proof of the purity theorem
2.2.1.

Theorem \ref{ffgg} is thus   proved in the case when $X$ is 
nonsingular. In order to tackle the case when $X$ is singular,
we must first endow intersection cohomology groups with a MHS, i.e. 
 we must now prove Theorem
\ref{tm111}.

\item
Proof of
Theorem \ref{tm111}. The projective case is proved in 
\ci{decmightam}. We follow the same strategy and
point out the needed modifications.
Let $f: X\to Y$ be  a resolution of the singularities of $Y$.
By the decomposition theorem, 
the intersection cohomology groups $I\!H^*(Y)$ and $I\!H^*_c(Y)$ are the subspaces
of the quotient perverse cohomology groups
 $H^*_0(X)$ and $H^*_c(X)$, respectively,  that correspond to the unique 
dense stratum on $Y$.  Part 1 follows by applying 
Theorem \ref{ffgg}, part 3, which we have proved in the case when 
$X$  is nonsingular.  
 The MHS so-obtained are  shown to be independent of 
 the resolution by an argument identical
 to the one in the proof of   \ci{decmightam}, proof of Theorem 2.2.3.a. 
The proof of part 2 follows from Theorem \ref{ffgg}, part 4 applied to $f$ ($X$ is nonsingular)
when we consider the dense stratum on $Y$. The proof of part 3 is identical
to the one for the projective case.

We can now  complete the proof of to Theorem \ref{ffgg} by dealing with 
the case when
$X$ is singular.

\item
Theorem \ref{ffgg}, part 1 follows directly from   Proposition
\ref{cor-555}. This proposition is the intersection cohomology analogue
of  Theorem \ref{hmhspl} and it is proved in pretty much the same way.
We only need to verify that
by taking general linear sections, the restriction maps on the  intersection cohomology
groups, 
and the  co-restriction maps on the intersection cohomology groups with compact 
supports are compatible with the 
MHS of Theorem \ref{tm111}. 

\item
Theorem \ref{ffgg}, part 1 clearly implies Theorem \ref{ffgg}, part 2.

\item 
Proof of Theorem \ref{ffgg}, part 3.
We discuss the case of cohomology. The case of cohomology with compact supports
is analogous.
Let $g: X' \to X$ be a resolution of the singularities of $X$. Set $h:= f\circ g: X' \to Y$.
Since $X'$ is nonsingular,  Theorem \ref{ffgg} holds for $g$ and for $h$.
Let $F^{g}_a$ and  $F^h_b$ denote the increasing perverse Leray filtrations
on the groups $H^*(X')$ associated with the maps $g$ and $h$, respectively.
Let $Gr^{F^g}_a H^*(X')$,  $Gr^{F^h}_b H^*(X')$ denote the
 corresponding graded pieces. 
 By Theorem \ref{ffgg}, part 1, the subspaces 
of both  filtrations are MHSS. 
 The graded pieces, as well
 as all  the bi-graded pieces  $Gr^{F^g}_a Gr^{F^h}_b H^*(X')$
 inherit the natural subquotient MHS.
 The same is true, via the just-established Theorem \ref{ffgg}, part 1, for the groups
 $I\!H^*(X)$,  for the subspaces
 $I\!H^*_{\leq i}(X)$  of the perverse Leray filtration
 with respect to $f$, and for their graded pieces
 $I\!H_i(X)$.
 Given  two finite  filtrations, $F$ and $G$  on an object $M$ of an abelian category,
Zassenhaus Lemma yields a canonical isomorphism $Gr^F_a Gr^G_b M \simeq 
Gr^G_b Gr^F_a M$. 
We apply this to the MHS $H^*(X')$ and obtain 
 a canonical isomorphism $Gr^{F^g}_a Gr^{F^h}_b H^*(X')
 \simeq
 Gr^{F^h}_b Gr^{F^g}_a H^*(X')$ of MHS.
 Note that the decomposition theorem
implies that each summand decomposes (as a vector space, a priori not
 as a MHS) according to the strata on $Y$ of a common refinement
of stratifications
of the maps $h$ and $f$, and that the Zassenhaus isomorphism
is a direct sum map (of vector spaces, a priori not of MHS).
 We have the canonical epimorphism of MHS
\[  Gr^{F^g}_0 H^*(X') \lorw  I\!H^* (X) \]
(this is how the MHS on the rhs has been constructed in the proof
of  Theorem \ref{tm111} given in n. 3 above).
This  map induces the epimorphic map of MHS
\[ 
Gr^{F^h}_i Gr^{F^g}_0 H^*(X') \lorw  I\!H_i^* (X).\]
This map is a direct sum map with respect to the strata $S$ on $Y$.
It remains to show that each $S$-summand on the left-hand-side is a MHSS.
Theorem \ref{ffgg}, part 3 (which we have proved above  for $X'$ is nonsingular)
implies  that $Gr^{F^h}_i  H^*(X')$ splits according to strata
into MHSS. Each $S$-summand of this group
maps onto the corresponding $S$-summand
in 
$Gr^{F^g}_0 Gr^{F^h}_i H^*(X')$ which is then a MHSS of this group.
We conclude by Zassenhaus Lemma, for the identification given
by  this lemma is compatible with the direct sum decomposition by strata
and with the MHS.

\item
Theorem \ref{ffgg}, part 4  is proved using the same argument employed 
in the case when $X$ is nonsingular (see item 1 of this list), provided we replace
$f_* \rat_X[n]$ with the self-dual  $f_* IC_X$ in Proposition \ref{pfbmhs}.

\item
The proof of 
Theorem \ref{appl3} in the case when $X$ and $Y$ are projective
is the main result of  \ci{decmigseattle}. The arguments provided
in that paper  are  quite general and work verbatim in the quasi projective case
and also in the case of compact supports (one only has to replace the 
pure Hodge structures employed there, with the MHS introduced here). 

\item

Intersection cohomology is not functorial in the 
``space'' variable. The paper \ci{cinque} constructs,
for every proper map $f: X \to Y,$
a non canonical map $I\!H^* (Y) \to I\!H^*(X).$
If $f$ is surjective, these morphisms stem from the decomposition theorem 
and are splitting injections. 
The same proof as \ci{decmigseattle}, Theorem 3.4.1, again replacing pure with mixed,
yields the proof of Theorem \ref{appl5}. 
\end{enumerate}

It is now easy to prove Theorem \ref{coincide}, i.e. to prove that
the mixed Hodge structures we construct coincide with the ones arising from M. Saito
theory of mixed Hodge modules (MHM).

\medskip
\n
{\bf Proof of Theorem \ref{coincide}.}

{\em First proof.}
Let $Y$ be a quasi projective irreducible variety. 
By \ci{saitomathann},
the MHS on
$H^*(Y)$ and $H^*_c(Y)$ stemming from the theory of MHM coincides
with the  one constructed by Deligne in \ci{ho3}. 
It follows that if
$X$ is nonsingular, then the two kinds of MHS
of  the subspaces  appearing in Theorem \ref{ffgg}, part 3 
coincide.   We apply this fact to a resolution of the singularities 
$f: X \to Y$ of $Y$ and we see  that the  two possible
MHS on the intersection cohomology groups $I\!H^*(Y)$ and $I\!H^*_c(Y)$
of an irreducible quasi projective variety coincide. Once the two MHS coincide
on the intersection cohomology groups, they also coincide 
on the subspaces appearing in Theorem \ref{ffgg} (and now $X$ can be singular).

{\em Second proof.} The splittings of Theorem \ref{appl5} arrive 
from a construction in homological algebra (due to Deligne; see \ci{decmigseattle})
that works in the category ${\mathcal D}_Y$ as well as in the derived category
$D^b(MHM_Y)$  of MHM
on $Y$. This means that we can take the splitting $\phi$
in ${\mathcal D}_Y$ to be the trace of a splitting in $D^b(MHM_Y)$,
and  this implies the conclusion.
\blacksquare

\medskip
This ends the outline of the proofs.  The next section is  devoted to 
completing the proofs.

Let me try to give an idea of how the mixed  Hodge-theoretic results are proved
 by looking at the  following special  simple case. Let
 $f: X \to Y$ be a  proper birational map of surfaces, with $X$ nonsingular, 
 such that $f$ is an isomorphism away from a 
 curve   $E\subseteq X$ contracted by $f$ to a point
 on $Y$.   The decomposition theorem implies that
 \[H^2(X)= I\!H^2(Y) \oplus \langle[E]\rangle, 
 \qquad 
 H^2_c(X)= I\!H_c^2(Y) \oplus \langle[E]\rangle.
 \]
 where $[E]$ is the fundamental class of $E$. 
 
 Let me illustrate the technique used in this paper to  
verify that in both equations  both summands are MHSS 
 of the usual MHS.

Poincar\'e duality yields an isomorphism
$ \iota: H^2(X) \simeq H^2_c(X)^{\vee}$.   The map $\iota$ is a direct sum map
with respect to the direct sum decompositions above.
 
 By invoking the decomposition theorem (see \ci{decmightam}, proof of the 
purity theorem 2.2.1, especially (44)),
  the pull-back map in cohomology 
 $r: H^2(X) \to H^2(E)$ is injective when restricted to  the
 summand $\langle [E] \rangle$ and it is the zero map
 on $I\!H^2(Y)$.
 Since $r$ is   a map of  MHS, its  kernel $I\!H^2(Y)$ is a MHSS of the MHS $H^2(X)$.
 
 We argue in the same way for cohomology with compact supports, except
 that the map $r$ is now the map $r': H^2_c(X) \to H^2_c(E)$
 which is the dual of  the proper push-forward map 
 $H_2^{BM}(E) \to H^{BM}_2(X)$ in Borel-Moore homology.
It follows  that $I\!H^2_c(Y)$ is a MHSS of  the MHS $H^2_c(X)$.
 
 The Poincar\'e isomorphism is in fact an isomorphism,
 $ \iota: H^2(X) \simeq H^2_c(X)^{\vee}(-2)$ of MHS.
 
The composition of maps of MHS
 $H^2(X) \to H^2_c(X)^{\vee}(-2) \to I\!H^2_c(Y)(-2)$ has kernel
 $\langle [E] \rangle$ which is then a MHSS of $H^2(X)$. Similarly, we show that
 $\langle [E] \rangle$  is  a MHSS of $H^2_c(X)$.

\subsection{Completion of the proofs of the results of $\S$\ref{dtamhs}}
\label{-1}
\begin{pr}
\label{pfbmhs}
Let $f:X \to Y$ be a proper map of 
irreducible quasi projective   varieties. Assume that  $X$  is nonsingular.
For every $i$ and $j$:

\begin{enumerate}
\item
the vector spaces  $H^j_i (X)$ and
$H^j_{c, i} (X)$ carry  natural MHS
which are  subquotients of the canonical MHS $H^*(X,
\rat)$ and $H^*_c(X, \rat)$. 

\item 
the Poincar\'e Pairing
$H^j (X) \simeq  H^{2m-j}_c  (X)^{\vee}(-m)$ descends 
to an   isomorphism of MHS
\[ H^j_{i}(X) \simeq  H^{2m-j}_{c,-i}  (X)^{\vee} (-m),
\qquad
H^j_{i,l,S}(X) \simeq  H^{2m-j}_{c,-i,l,S}  (X)^{\vee} (-m)
.\]
\end{enumerate}
\end{pr}
{\em Proof.}
The spaces in question are the graded pieces of the perverse Leray filtration
and part 1  follows from Theorem \ref{hmhspl}. 

We turn to part 2. By  the mixed Hodge theory of algebraic varieties (\ci{ho3}), 
the Poincar\'e pairings  $H^j(X) \simeq H^{2m-j}_c (X)^{\vee} (-m)$
is an isomorphism of MHS.
By  \ci{decmightamv1}, Lemma 2.9.1
(which proof, written for $X$ proper and nonsingular, is valid when $X$ is merely nonsingular), the Poincar\'e pairing above is compatible with the perverse Leray filtration,
i.e. for $i\in \zed$ it induces   maps $H^j_{\leq i}(X) \to 
H^{2m-j}_{c, \leq -i}(X)^{\vee}(-m)$, 
and  it 
 descends to each $i$-th graded group as a  linear isomorphisms
 $P: H^j_i(X) \simeq H^{2m-j}_{c,-i}(X)^{\vee}(-m)$. This linear isomorphism
 is of MHS, for the subquotient MHS on the graded groups.
 This establishes the first statement of part 2, i.e.  for $H^j_i(X)$.  
 
 We now turn to  $H^j_{i,l,S}(X)$.
 By the same lemma quoted above,
the linear isomorphism $P$  coincides with the map in hypercohomology associated with  
 the canonical isomorphism stemming from Verdier Duality (recall that $f_*\rat_X[m]$
is self-dual)
\[
\pc{i}{f_*\rat_X[m]} \simeq \pc{-i}{f_* \rat_{X[m]}}^{\vee}.
\]
Both  sides split according to strata as in (\ref{grouper}).
We are left with showing that the map $P$ is a direct sum map for
this decomposition according to strata. This follows immediately from the fact that there are no
nontrivial maps between intersection complexes supported on different
subvarieties.
\blacksquare

\bigskip
The proof    of  the two lemmata below is a combination of the 
results of this paper and of the
methods employed in \ci{decmightam} 
in  the proof the purity Theorem 2.2.1. The key new ingredient  is 
Proposition \ref{pfbmhs}. Another difference is that, since we need to argue
using pairings, we need to simultaneously  
keep track of cohomology and of cohomology
with compact supports, even if we are interested only in cohomology.

Even though we do not repeat the parts of the proof that are contained in \ci{decmightam},
for the reader's convenience in the course of the proof we quote the relevant
results from \ci{decmightam}

As in the proof of the purity Theorem 2.2.1 in \ci{decmightam},
the proof is by induction $m:=\dim{X}$.    The first main   step  is carried out
in the following lemma,  where we deal with the cases $(i,j) \neq (0,m)$
as follows: i) we  take hyperplane sections of the domain
$X$ and   deal with the cases when $i\neq 0$,
and ii) we take
hyperplane sections of the target
$Y$ and  deal with the cases $(i=0, j\neq m)$.
The remaining and crucial case when $(i=0,j=m)$  is dealt-with in the second main step,
Lemma \ref{crucial}

\begin{lm}
\label{indst}
If Theorem {\rm \ref{ffgg}}, part {\rm 3}  holds for every projective map
$g: Z \to Z'$ of quasi-projective varieties, $Z$ nonsingular, $\dim{Z} < \dim {X},$
then it holds for every group  $H^j_i(X)$ and $H^j_{c,i}(X)$ with
$(i,j) \neq (0,m).$
\end{lm}
{\em Proof.}  
See \ci{decmightam}, Lemma 7.1.1. Choose a general hyperplane section $r : X^1 \to X.$
In this proof we need the section to be smooth, so that
we can  apply induction, and transverse
 to the relevant stratifications,  so that $r^! \simeq r^*[-2]$  when applied to
 complexes constructible with respect to those same stratifications.
In what follows, we shall use freely these 
facts as well as that  $f_{*} =f_!,$ and $g_{*}= g_!.$
We have  the proper  map  $g: X^1 \to Y $ and
the affine map $u: X \setminus X^1 \to Y$.
There is the adjunction map \[f_{*} \rat_X[m] \lorw g_{*} \rat_{X^1}[m-1] [1].\]
Taking hypercohomology and hypercohomology with compact supports, we obtain 
  restriction maps \[H^j_i(X) \lorw H^j_{i+1} (X^1), \qquad  
H^j_{c,i}(X) \to H^j_{c , i+1} (X^1).\] Similarly,  by considering the adjunction map
\[g_{*} \rat_{X^1}[m-1][-1] \lorw f_{*} \rat_X[m],\] we obtain Gysin maps
\[H^{j-2}_{i-1}(X^1) \lorw H^j_{i}(X), \qquad 
H^{j-2}_{c,i-1}(X^1) \lorw H^j_{c } (X).\] 
By the Weak-Lefschetz-type Proposition 4.7.6 in \ci{decmightam}
(the key point is that $u$ is affine, hence $t$-left exact), we have
that the natural restriction-type map
$$
\pc{i}{ (f_{*} \rat_X[m] )} \lorw  \, \pc{i+1}{ (g_{*} \rat_{X^1} [m-1] )}
$$
is a splitting monomorphism  for $i <0,$ and that the natural 
Gysin-type map
$$
\pc{i-1}{ (g_{*} \rat_{X^1} [m-1] ) } \lorw    \, \pc{i}{ (f_{*} \rat_{X} [m] ) }
$$
is a splitting epimorphism for $i>0$.

Since the restriction and Gysin maps above are direct sum maps with respect
to the direct sum decomposition by strata, 
the statement of the  lemma  follows  for every $(i,j)$ with $i \neq 0$ by virtue of the
 the hypotheses applied to $Z=X^1 \to Y = Z'$ (cfr. \ci{decmightam}, p.744, bottom).
 
\smallskip
Let $i=0$ and $ j \neq m.$ The argument is formally similar to the one just given.
However, instead of using   hyperplane sections of $X$ and 
Weak-Lefschetz-type results for 
the  perverse cohomology complexes $\pc{i}{ (f_{*} \rat_X[m] )}$,
$i \neq 0$, we work with 
hyperplane sections on $Y$ and Weak-Lefschetz-type results
on the cohomology  groups \[H^j(Y, \pc{0}{ (f_{*} \rat_X[m] )}, \qquad
H^j_c(Y, \pc{0}{ (f_{*} \rat_X[m] )}, \qquad \forall j \neq 0.  \]

We need to show that
 the theorem holds for the following four groups.
 \[ 1) \;H^{j<m}_0(X), \quad   2) \; H^{j>m}_{c,0}(X), \quad
  3) \; H^{j>m}_0(X),  \quad  4)  \; H^{j<m}_{c,0}(X) .\]

  The cases 1) and 2) are dual to each other and so are 3) and 4). 
  It follows that it is enough to establish the result
  in cases 1) and 3).

  We choose  an affine embedding $Y$ into projective space
  and a general hyperplane
 section $Y_1$ of $Y$  with respect to this embedding. Note that, in this case,
 $Y\setminus Y_1$ is affine.
We have the
 associated map $h:X_1:=f^{-1}(Y_1) \to Y_1$ and the decomposition theorem
 for $h$ takes the form of the decomposition  (\ref{dtgl}) restricted to $Y_1$
 and shifted by $[-1]$ and
 we have 
 \[\pc{i}{ h_{*} \rat_{X_1}[m-1]}= \pc{i}{ f_{*} \rat_X[m]}_{|Y_1} [-1].\]
 Note that the skyscraper summands disappear
 after restriction (though this plays no role in the sequel of this proof).

\smallskip 
  Case 1).  
 The natural restriction maps  $H^j(X) \to H^j(X_1)$ are of MHS.
 It follows that the maps 
   $H^j_0(X) \to H^j_0(X_1)$ induced 
  on the subquotients  are of MHS with respect to the MHS stemming from Proposition
  \ref{pfbmhs}. Moreover, 
  these induced maps
 are direct sum maps with respect to the direct sum decompositions by strata
 (\ref{mnb}).
  By  the inductive hypothesis, the theorem holds for $h: X_1 \to Y_1$. It is
  thus  enough 
 to show that the natural restriction maps above are injective  for $j <m$.
 Since  the Lefschetz hyperplane  Theorem \ref{swl} applied
 to the perverse sheaf $\pc{0}{ (f_{*} \rat_X[m]}$ on $Y$ implies injectivity in the desired range, the result follows in case 1).
 
 \smallskip
  Case 3). The  restriction map (see Remark \ref{wswl})
 \[H^j(X_1, i^! \rat)= H^{j-2}(X_1, \rat) \lorw H^j(X, \rat)\]
 is the  natural Gysin map and is of MHS.
 Passing to graded groups, the induced Gysin-type map $H^{j-2}_0 (X_1) \to H^j_0 (X)$
 is a map of MHS and it is a  direct sum map with respect to the decomposition by strata. As in  case 1), it is enough to show that
 this Gysin-type map is  surjective for $j>m$. The Gysin-type map in question  appears
 in  the long exact sequence of  cohomology of the triangle
 \[I_! I^!\pc{0}{ (f_{*} \rat_X[m]} \lorw \pc{0}{ (f_{*} \rat_X[m]}
 \lorw J_* J^* \pc{0}{ (f_{*} \rat_X[m]} \stackrel{+}\lorw  
 \]
  on $Y$, where $I: Y_1 \to Y \leftarrow Y-Y_1: J$ are the natural inclusions.
 To establish the wanted surjectivity, it  is enough to observe that
 \[H^r(Y, J_* J^* \pc{0}{ (f_{*} \rat_X[m]})=
 H^r(Y-Y_1, J^* \pc{0}{ (f_{*} \rat_X[m]}) =0 \]
 by the theorem on the cohomological dimension of affine varieties
 for  perverse sheaves applied to the perverse  $J^* \pc{0}{ (f_{*} \rat_X[m]}$.
 The result follows also  in case 3).
  \blacksquare

\medskip
The  following lemma takes care of   the remaining cases $H^m_0(X)$ and $H^{m}_{c,0}(X)$
and completes the proof of Theorem \ref{ffgg}, part 3, in the case when $X$ is nonsingular.

\begin{lm}\label{crucial}
Theorem {\rm \ref{ffgg}}, part {\rm 3} holds when $X$ is nonsingular.
\end{lm}
{\em Proof.} The proof is by induction on $\dim{X}$.
The cases $\dim{X=0, 1}$ are trivial.

We assume that the theorem holds for every map $g: Z \to Z'$ as in  Lemma \ref{indst}.
By the same lemma, we are left with the cases of $H^m_0(X)$ and $H^m_{c,0}(X).$

We choose a non-dense stratum  $S$ in $f(X)$  and 
we proceed exactly as in the  proof of Theorem 2.2.1 in \ci{decmightam}  and we prove that
\begin{equation}
\label{xxx}
\bigoplus_{l, S' \neq S} H^m_{0,l,S'}(X) \subseteq H^m_0(X)
\end{equation}
is a MHSS.  The only minor difference is that when we take
the closure $Z'_S$  of $f^{-1}(S)$ in $X$ and  take a resolution
of the singularities  $\rho: Z_S \to Z'_S$, the resulting
quasi projective variety $Z_S$ is not projective, however, the properness
of $Z_S$ plays no role in \ci{decmightam}, what is essential
is the fact that  $Z'_S \to  \overline{S}$ is proper. 
Let us point out that  we can  use the inductive hypothesis in view of the fact that
$\dim{Z_S} < \dim{X}$,
and this explains why we started with a  non-dense stratum.

By taking intersections, we see  that  any  direct sum of terms
which includes the  dense stratum $\Sigma \subseteq f(X)$
gives  a MHSS of $ H^m_0(X).$

The same line of reasoning works for $H^m_{c,0}(X)$. The only 
difference is 
that the maps we use are  not a pull-back maps in cohomology, but 
rather the maps in cohomology with compact supports  which are the  duals
of the proper push-forward maps  in Borel-Moore homology.


We are left with the case of the summands associated with the dense
stratum $\Sigma$ in $f(X)$.
Consider the composition of maps of MHS (dualizing turns a MHSS into a quotient MHS):
$$
H^m_0 (X) \simeq H^m_{c,0}(X)^{\vee} (-m) \lorw  \bigoplus_{l, S' \neq \Sigma}{
H^m_{c,0, l,S'} (X)^{\vee} (-m)}.
$$
The kernel, $H^m_{0,\dim{\Sigma},\Sigma}(X)$, is a MHSS of $H^m_0(X)$.
It follows that all direct summands $H^m_{0,l,S}(X),$
of $H^m_0(X)$ are MHSS. By dualizing, 
the same  is true for $H^m_{c,0}(X)$ and its summands.
\blacksquare

\bigskip
The following is the intersection cohomology analogue of Theorem
\ref{hmhspl}. It is needed to endow the perverse cohomology  groups
$I\!H_i^*(X)$ and $I\!H^*_{c,i}(X)$ with MHS.
 For convenience only, in the proof below,
we use a different normalization for the intersection
complex, i.e. if  $X$ is  an irreducible variety, then set ${\mathcal IC}_X
 := IC_X [-\dim{X}]$ so that
 $I\!H^j(X)= H^j(X, {\mathcal IC}_X)$.  If $X$ is  irreducible and nonsingular, then
 ${\mathcal IC}_X = \rat_X$.

\begin{pr}
\label{cor-555}
{\bf (Perverse Leray and MHS: intersection cohomology)}
Let $f:X \to Y$ be a map of irreducible quasi projective    varieties.
There are an integer $d$, a  variety ${\mathcal X}$ and   a flag ${\mathcal X}_*$
on it
 such that there are identities of   filtered groups
$$
(I\!H^*(X), L_{\pp}^f  ) \; = \; ( I\!H^*({\mathcal X}), F_{ {\mathcal X}_* } (-d)),
\qquad
(I\!H^*_c(X), L_{\pp}^f  ) \; = \; (I\!H^*_c({\mathcal X}), G_{ {\mathcal X}_* }(d) ).
$$
In particular,  the subspaces of the perverse Leray filtrations on $I\!H^*(X)$ and $I\!H^*_c(X)$ are  MHSS
(for the MHS of Theorem \ref{tm111}), and the graded groups
 $I\!H_i(X)$ and $I\!H_{c,i}(X)$ inherit the subquotient  MHS.
\end{pr}
{\em Proof.} We employ  the set-up of $\S$\ref{ssjo}, especially
(\ref{cdjo}). Note that $\pi^* {\mathcal IC}_X = {\mathcal IC}_{\mathcal X}$.
The identities of filtered groups stem from Theorem  \ref{flplr}
applied to $C= {\mathcal IC}_X$.

CLAIM: Let $i:  {\mathcal H} \to {\mathcal Y}$ be 
a general codimension $c$  linear section (for any embedding in projective space).
Let $i: {\mathcal X}_H \to {\mathcal H}$ be  the pre-image of $\mathcal H$.
The restriction  map $I\!H^*( {\mathcal X}) \to I\!H^*({\mathcal X}_{\mathcal H})$
and the co-restriction map $I\!H^*_c( {\mathcal X}_{\mathcal H}) \to I\!H^*_c({\mathcal X})$ 
are maps of MHS for the MHS of Theorem \ref{tm111}.

{\em Proof of the {\rm CLAIM}}.  
 Let $h: {\mathcal X}' \to {\mathcal X}$ be a resolution of the singularities
of ${\mathcal X}$. 
We have the  Cartesian diagram
\[
\xymatrix{
{\mathcal X}'_{\mathcal H}   \ar@/^2pc/[rr]^{f':= f\circ h}
 \ar[r]^h \ar[d]^i & {\mathcal X}_{\mathcal H} \ar[r]^f \ar[d]^i & {\mathcal H}  \ar[d]^i \\
{\mathcal X}' \ar[r]^h \ar@/_2pc/[rr]^{f':= f\circ h}  & {\mathcal X} \ar[r]^f  & {\mathcal Y}. 
}
\]
For ${\mathcal H}$ sufficiently general, 
${\mathcal IC}_{{\mathcal X}_{\mathcal H}} = i^* {\mathcal IC}_{\mathcal X}$
and $f': {\mathcal X}'_{\mathcal H} \to {\mathcal X}_{\mathcal H} $ is also  a resolution.

By the decomposition theorem,
the intersection complex  ${\mathcal IC}_{{\mathcal X}}$ is a direct summand of
$h_* \rat_{{\mathcal X}'}$ (see \ref{dtgl}).   Similarly, for ${\mathcal IC}_{{\mathcal X}_{\mathcal H}}$
and $h_* \rat_{ {\mathcal X}'_{\mathcal H}   }$.
The two direct sum decompositions are related by the restriction map
$h_* \rat_{{\mathcal X}'} \to i_*h_* \rat_{{\mathcal X}'_{\mathcal H}}$. This map is a direct sum 
map and maps ${\mathcal IC}_{{\mathcal X}}$ to ${\mathcal IC}_{ {\mathcal X}_{\mathcal H}}$.
This map yields the restriction map in intersection cohomology
$I\!H^*({\mathcal X}) \to I\!H^*({\mathcal X}_{\mathcal H})$.
Similarly, for the co-restriciton map.

We complete the proof of the CLAIM for the restriction map. The case of the co-restriction map 
is completely analogous. The restriction map  $H^*({\mathcal X}' ) \to H^*({\mathcal X}'_{\mathcal H})$
is of MHS. The spaces  $H_{\leq 0}({\mathcal X} ) \subseteq H^*({\mathcal X} )$ and
$H_{\leq 0}({\mathcal X}'_{\mathcal H} ) \subseteq H_{\leq 0}({\mathcal X}'_{\mathcal H} )$
are MHSS, mapped into each other via the restriction map.
It follows that the restriction map descends to a map
 $H_0({\mathcal X}') \to H_0({\mathcal X}'_{\mathcal H} )$ of MHS. 
 This map is a direct sum map with respect to the strata and, by Theorem
\ref{ffgg}, part 3,   each component 
is a map of MHS. Since the restriction map in intersection cohomology
is one of these summands, the CLAIM follows.

We conclude the proof by observing that the flag ${\mathcal X}_*$ is a pull-back flag
of a general flag ${\mathcal Y}_*$ on ${\mathcal Y}$ and, by applying the CLAIM
to the elements of the flag, we see that the kernels (images, resp.)
of the restriction (co-restriction, resp.)  maps, i.e. the subspaces of the perverse filtrations,
  are MHSS.
\blacksquare

\subsection{The results of $\S$\ref{dtamhs} hold for  not necessarily
 irreducible varieties}
\label{niv}
In this section we point out that, while we have stated most of the Hodge-theoretic applications
in $\S$\ref{dtamhs} in the context of  quasi projective {\em irreducible} varieties, these results in fact hold
for quasi projective varieties. 

The key point is to give the correct definition of intersection complex
for not necessarily irreducible, nor pure dimensional varieties. Once that definition
is in place, the rest follows quite easily.

Given an irreducible variety $Y$ of dimension $n$, one defines the intersection complex
$IC_Y$ as follows: let $j:U \subseteq Y$ be a Zariski-dense, open
and nonsingular subset of $Y$; then $IC_Y = j_{!*} \rat_U [n]$ is the intermediate
extension of the perverse sheaf  $\rat_U [n]$.

Let $Y$ be any variety and $Y_{reg} = \coprod_{l\geq 0} U_l$
be the decomposition of the regular part in 
pure $l$-dimensional components. 
Let ${\frak Q} := \bigoplus_{l\geq 0} \rat_{U_l}[l]$ and $j: Y_{reg} \to Y$ be the open immersion.
 Define the intersection complex of $Y$ as follows
 \[IC_Y =  j_{!*} {\frak Q}.\]
 
 The following are easily verified:
 \begin{enumerate}
 \item
 $IC_Y = \bigoplus_{l\geq 0} IC_{\overline{U}_l}$;
 \item
 let $\nu: Y' \to Y$ be the normalization; recall
 that $Y'$ is a disjoint union of irreducible normal varieties and that $\nu$ is finite,
 so that $\nu_*= R^0 \nu_*$;
  we have
 $\nu_* IC_{Y'} = IC_{Y}$. 
 \end{enumerate}
 
 It follows that with this definition of intersection complex
  the decomposition theorem holds for a proper map of
  varieties. In fact, one normalizes the domain and works on each component separately.
  
  It is easy to generalize all the mixed-Hodge-theoretic applications
  in $\S$\ref{dtamhs}
  of this paper to arbitrary quasi projective varieties, provided we use
  the definition of intersection complex given above.
  
  We leave to the reader the task of verifying that 
  all the  proofs go through verbatim, with the possible exception
  of the one of Theorem \ref{tm111}, part 3.  
  
  In this case,
  we need to verify that
  the map $a:H^*(Y) \to I\!H^*(Y)$ is of MHS. This is done by reduction
  to the irreducible case as follows. Let $\{Y_t\}$ be the set of irreducible components of $Y$.
  We have that $I\!H^*(Y)= \oplus_t I\!H^*(Y_t) $.  The map $a$ is of MHS iff
  the induced maps $a_t: H^*(Y) \to I\!H^*(Y_t)$  are all of MHS.
  The map $a_t$ factors as follows $H^*(Y) \to H^*(Y_t) \to I\!H^*(Y_t)$. The first map
  if of MHS by the theory of MHS for cohomology. The second one is
  of MHS by Theorem \ref{tm111}, part 3.

\bibliographystyle{amsplain}

\end{document}